\documentclass[10pt]{article}
\usepackage{graphicx}
\usepackage{amsmath}
\usepackage{amsfonts}
\usepackage{amssymb}
\begin{document}

\title{Polynomial decay rate \\for the dissipative wave equation}
\author{Kim Dang \textsc{Phung}\\{\small 17 rue L\'{e}onard Mafrand 92320 Ch\^{a}tillon, France}\qquad\\phung@cmla.ens-cachan.fr}
\date{}
\maketitle

\bigskip

\section{Introduction and main result}

\bigskip

This paper is devoted to study the stabilization of the linear wave equation
in a bounded domain damped in a subdomain when the geometrical control
condition (see \cite{BLR}) of the work of C. Bardos, G. Lebeau and J. Rauch is
not fulfilled. In such case, they \cite{BLR} proved that the uniform
exponential decay rate of the energy cannot be hoped due to the existence of a
trapped ray that never reaches the support of the damping. Another important
contribution in this field was done by G. Lebeau \cite{Le} who establishes a
logarithmic decay rate for the dissipative wave equation when no assumption on
rays of geometrical optics is required, but when more regularity on the
initial data is allowed. Now, it seems natural to search a general description
of the geometries of both domain and support of the damping under which the
energy of the dissipative wave equation decreases in a polynomial way. A first
answer in this direction was done by Z. Liu and B. Rao \cite{LR} who consider
the wave equation on a square damped in a vertical strip. In this paper, we
improve the geometry to a partially cubic domain where the damping acts in a
neighborhood of the boundary except between a pair of parallel square face of
the cube.

\bigskip

Before stating the main result of this paper, we begin by presenting precisely
the geometry of our problem. Next, we introduce the equations that will be
used throughout this work. Our main result is given at the end of this section.

\bigskip

\subsection{The geometry}

\bigskip

Let $m_{1},m_{2},\rho>0$. Let $\Omega$ be a connected domain in $\mathbb{R}%
^{3}$ bounded by $\Gamma_{1}$, $\Gamma_{2}$, $\Upsilon$ where

\begin{description}
\item $\quad\Gamma_{1}=\left[  -m_{1},m_{1}\right]  \times\left[  -m_{2}%
,m_{2}\right]  \times\left\{  \rho\right\}  $, with boundary $\partial
\Gamma_{1}$,

\item $\quad\Gamma_{2}=\left[  -m_{1},m_{1}\right]  \times\left[  -m_{2}%
,m_{2}\right]  \times\left\{  -\rho\right\}  $, with boundary $\partial
\Gamma_{2}$,

\item $\quad\Upsilon$ is a surface such that $\Upsilon\subset\mathbb{R}%
^{2}\left\backslash \left(  \left(  -m_{1},m_{1}\right)  \times\left(
-m_{2},m_{2}\right)  \right)  \right.  \times\mathbb{R}$, with boundary
$\partial\Upsilon=\partial\Gamma_{1}\cup\partial\Gamma_{2}$.
\end{description}

\noindent We suppose that the boundary of $\Omega$, $\partial\Omega=\Gamma
_{1}\cup\Gamma_{2}\cup\Upsilon$, is $C^{\infty}$. Let $\Theta$ be a
neighborhood of $\Upsilon$ in $\mathbb{R}^{3}$ and let $\omega=\Omega
\cap\Theta$.

\bigskip

\noindent As $\Theta$ is a neighborhood of $\Upsilon$ in $\mathbb{R}^{3}$,
there exists $r_{o}\in\left(  0,\min\left(  m_{1},m_{2},\rho\right)
/2\right)  $ such that $\left[  x_{b1}-2r_{o},x_{b1}+2r_{o}\right]
\times\left[  x_{b2}-2r_{o},x_{b2}+2r_{o}\right]  \times\left[  x_{b3}%
-2r_{o},x_{b3}+2r_{o}\right]  \subset\Theta$ for any $\left(  x_{b1}%
,x_{b2},x_{b3}\right)  \in\partial\Gamma_{1}$. Let $h_{o}=\min\left(
1,\left(  r_{o}/8\right)  ^{2}\right)  $.

\noindent Next, we choose $\omega_{o}=\left(  -m_{1}+r_{o},m_{1}-r_{o}\right)
\times\left(  -m_{2}+r_{o},m_{2}-r_{o}\right)  \times\left(  -\frac{\rho}%
{4},\frac{\rho}{4}\right)  $.

\bigskip

Throughout this paper, $c$ denotes a positive constant which only may depend
on $\left(  m_{1},m_{2},\rho\right)  $. Also $\gamma$ will denote an absolute
constant larger than one. The value of $c>0$ and $\gamma>1$ may change from
line to line.

\bigskip

\subsection{The equations}

\bigskip

We consider the dissipative wave equation in $\Omega$ with initial data
$\left(  w_{0},w_{1}\right)  \in H^{2}\left(  \Omega\right)  \cap H_{0}%
^{1}\left(  \Omega\right)  \times H_{0}^{1}\left(  \Omega\right)  $.%
\begin{equation}
\left\{
\begin{array}
[c]{rll}%
\partial_{t}^{2}w-\Delta w+\alpha\left(  x\right)  \partial_{t}w & =0 &
\quad\text{in}~\Omega\times\mathbb{R}_{+}\\
w & =0 & \quad\text{on}~\partial\Omega\times\mathbb{R}_{+}\\
\left(  w\left(  \cdot,0\right)  ,\partial_{t}w\left(  \cdot,0\right)  \right)
& =\left(  w_{0},w_{1}\right)  & \quad\text{in}~\Omega\text{ ,}%
\end{array}
\right.  \tag{1.1}\label{1.1}%
\end{equation}

\noindent with a non-negative dissipative potential $\alpha\in L^{\infty
}\left(  \Omega\right)  $ such that $\alpha>0$ in $\omega$. Denote%
\[
\mathcal{E}\left(  w,t\right)  =\int_{\Omega}\left(  \left|  \partial
_{t}w\left(  x,t\right)  \right|  ^{2}+\left|  \nabla w\left(  x,t\right)
\right|  ^{2}\right)  dx=\mathcal{E}\left(  w,0\right)  -2\int_{0}^{t}%
\int_{\Omega}\alpha\left(  x\right)  \left|  \partial_{t}w\left(
x,\theta\right)  \right|  ^{2}dxd\theta\text{ ,}%
\]

\noindent and recall that $\mathcal{E}\left(  w,t\right)  $ is a continuous
decreasing function of time.

\bigskip

\noindent Denote by $\{\mu_{j}\}_{j\geq1}$, $0<\mu_{1}<\mu_{2}\leq\mu_{3}%
\leq\cdot\cdot\cdot$, the eigenvalues of $-\Delta$ on $\Omega$ with Dirichlet
boundary conditions and by $\{\ell_{j}\}_{j\geq1}$ the corresponding
normalized eigenfunctions, i.e., $\Vert\ell_{j}\Vert_{L^{2}(\Omega)}=1$. Let
$u=u\left(  x,t\right)  $ be the solution of the wave equation
\begin{equation}
\left\{
\begin{array}
[c]{rll}%
\partial_{t}^{2}u-\Delta u & =0 & \quad\text{in}~\Omega\times\mathbb{R}\\
u & =0 & \quad\text{on}~\partial\Omega\times\mathbb{R}\\
\left(  u\left(  \cdot,0\right)  ,\partial_{t}u\left(  \cdot,0\right)  \right)
& =\left(  u_{0},u_{1}\right)  & \quad\text{in}~\Omega\text{ .}%
\end{array}
\right.  \tag{1.2}\label{1.2}%
\end{equation}

\noindent Suppose that $u_{0}=\sum\limits_{j\geq1}b_{j}^{0}\ell_{j}$ and
$u_{1}=\sum\limits_{j\geq1}b_{j}^{1}\ell_{j}$ are such that
\[
\left\|  u_{0}\right\|  _{H^{2}\cap H_{0}^{1}\left(  \Omega\right)  }^{2}%
=\sum\limits_{j\geq1}\mu_{j}^{2}\left|  b_{j}^{0}\right|  ^{2}<+\infty\text{
and }\left\|  u_{1}\right\|  _{H_{0}^{1}\left(  \Omega\right)  }^{2}%
=\sum\limits_{j\geq1}\mu_{j}\left|  b_{j}^{1}\right|  ^{2}<+\infty\text{ ,}%
\]

\noindent then it is known that
\[
u\left(  \cdot,t\right)  =\sum\limits_{j\geq1}\left[  b_{j}^{0}\cos\left(
t\sqrt{\mu_{j}}\right)  +\frac{b_{j}^{1}}{\sqrt{\mu_{j}}}\sin\left(
t\sqrt{\mu_{j}}\right)  \right]  ~\ell_{j}\text{ ,}%
\]

\noindent and $u\in C\left(  \mathbb{R};H^{2}\left(  \Omega\right)  \cap
H_{0}^{1}\left(  \Omega\right)  \right)  \cap C^{1}\left(  \mathbb{R}%
;H_{0}^{1}\left(  \Omega\right)  \right)  \cap C^{2}\left(  \mathbb{R}%
;L^{2}\left(  \Omega\right)  \right)  $. Let
\[
\mathcal{G}\left(  u,t\right)  =\int_{\Omega}\left(  \left|  \partial
_{t}u\left(  x,t\right)  \right|  ^{2}+\left|  \nabla u\left(  x,t\right)
\right|  ^{2}\right)  dx=\mathcal{G}\left(  u,0\right)  \text{ .}%
\]

\bigskip

\subsection{The bicharacteristics}

\bigskip

It is usual to associate with the wave equation, the geodesics. Recall that
the bicharacteristics are curves in the space-time variables and their Fourier
variables described by%
\[
\left\{
\begin{array}
[c]{rl}%
x\left(  s\right)  & =x_{o}+2\xi\left(  s\right)  s\\
t\left(  s\right)  & =t_{o}-2\tau\left(  s\right)  s
\end{array}
\right.  \text{\quad and }\left\{
\begin{array}
[c]{rl}%
\xi\left(  s\right)  & =\xi_{o}\\
\tau\left(  s\right)  & =\tau_{o}%
\end{array}
\right.
\]

\noindent with $\left|  \xi\left(  s\right)  \right|  ^{2}-\tau^{2}\left(
s\right)  =0$ for $s\in\left[  0,+\infty\right)  $, when $\left(  x_{o}%
,t_{o},\xi_{o},\tau_{o}\right)  \in\mathbb{R}^{4}\mathbb{\times R}%
^{4}\left\backslash \left\{  0\right\}  \right.  $. The rays are the
projection of the bicharacteristics on the space-time domain%
\[
\left\{
\begin{array}
[c]{rl}%
x\left(  s\right)  -x_{o}-2\xi_{o}s & =0\\
t\left(  s\right)  +2\tau_{o}s & =0\\
\left|  \xi_{o}\right|  ^{2}-\tau_{o}^{2} & =0\text{ ,}%
\end{array}
\right.
\]

\noindent here, $t_{o}=0$ and $\tau_{o}\neq0$. The generalized rays are rays
taking into account the geometry by following the rules of optic geometric.

\bigskip

\noindent The key geometric observation in our setting is that we do not have
that any generalized ray meets $\omega$. More precisely, any $\left(
x_{o},t_{o},\xi_{o},\tau_{o}\right)  \in\omega_{o}\times\mathbb{R\times R}%
^{4}\left\backslash \left\{  0\right\}  \right.  $ such that $\xi_{o}=\left(
0,0,\pm1\right)  $ generates a trapped generalized ray which never goes
outside $\left[  -m_{1},m_{1}\right]  \times\left[  -m_{2},m_{2}\right]
\times\left[  -\rho,\rho\right]  $. As a result, we do not have an uniform
exponential decay for the dissipative wave equation for any damping only
acting in $\omega$. Of course, the logarithmic decay rate still holds from
Carleman inequalities. However, we may hope a better decay rate because our
trapped generalized ray behaves quite simply by bouncing between $\Gamma_{1}$
and $\Gamma_{2}$ always in the same direction $\xi_{o}=\left(  0,0,\pm
1\right)  $. Furthermore, observe that the geometrical control condition is
fulfilled with $\omega\cup\omega_{o}$.

\bigskip

\subsection{Main result}

\bigskip

Now we are ready to formulate our main result.

\bigskip

Theorem.- \textit{There exist }$C>0$\textit{ and }$\delta>0$\textit{ such that
for any }$t>0$\textit{ and any initial data }$\left(  w_{0},w_{1}\right)  \in
H^{2}\left(  \Omega\right)  \cap H_{0}^{1}\left(  \Omega\right)  \times
H_{0}^{1}\left(  \Omega\right)  $\textit{, the solution }$w$\textit{ of
}(\ref{1.1})\textit{ satisfies,}
\begin{equation}
\int_{\Omega}\left(  \left|  \partial_{t}w\left(  x,t\right)  \right|
^{2}+\left|  \nabla w\left(  x,t\right)  \right|  ^{2}\right)  dx\leq\frac
{C}{t^{\delta}}\left\|  \left(  w_{0},w_{1}\right)  \right\|  _{H^{2}\cap
H_{0}^{1}\left(  \Omega\right)  \times H_{0}^{1}\left(  \Omega\right)  }%
^{2}\text{ .} \tag{1.3}\label{1.3}%
\end{equation}

\bigskip

Our strategy to get such polynomial decay rate will consist to establish an
observability estimate for the wave equation $u$ solution of (\ref{1.2}), or
more precisely, an inequality which traduces the unique continuation property
for $u$ between $\omega\times\left(  0,T\right)  $ and $\Omega\times\left\{
0\right\}  $ for some $T<+\infty$. For example, it is now known that an
interpolation inequality of the form
\[
\left\|  \left(  u_{0},u_{1}\right)  \right\|  _{H_{0}^{1}\left(
\Omega\right)  \times L^{2}\left(  \Omega\right)  }^{2}\leq C~\Lambda
^{1/\delta}\int_{0}^{T}\int_{\omega}\left|  \partial_{t}u\left(  x,t\right)
\right|  ^{2}dxdt
\]

\noindent where $\Lambda=\frac{\left\|  \left(  u_{0},u_{1}\right)  \right\|
_{H^{2}\cap H_{0}^{1}\left(  \Omega\right)  \times H_{0}^{1}\left(
\Omega\right)  }^{2}}{\left\|  \left(  u_{0},u_{1}\right)  \right\|
_{H_{0}^{1}\left(  \Omega\right)  \times L^{2}\left(  \Omega\right)  }^{2}}$,
implies (\ref{1.3}). On another hand, it is not difficult to deduce from
(\ref{1.3}) an inequality in the following form
\[
\left\|  \left(  u_{0},u_{1}\right)  \right\|  _{H_{0}^{1}\left(
\Omega\right)  \times L^{2}\left(  \Omega\right)  }^{2}\leq c\int_{0}^{\left(
C\Lambda\right)  ^{1/\delta}}\int_{\Omega}\alpha\left(  x\right)  \left|
\partial_{t}u\left(  x,t\right)  \right|  ^{2}dxdt
\]

\noindent and conversely, the later inequality implies (\ref{1.3}) where the
constant $C$ may change of values. In this paper, (\ref{1.3}) comes from
\[
\left\|  \left(  u_{0},u_{1}\right)  \right\|  _{H_{0}^{1}\left(
\Omega\right)  \times L^{2}\left(  \Omega\right)  }^{2}\leq C~\Lambda^{\gamma
}\int_{0}^{\left(  C\Lambda\right)  ^{\gamma}}\int_{\omega}\left|
\partial_{t}u\left(  x,t\right)  \right|  ^{2}dxdt\text{ .}%
\]

\bigskip

As it was suggested in the previous subsection, we have to pay more attention
on a ray of geometrical optics bouncing up and down infinitely between the two
parallel planes $\Gamma_{1}$ and $\Gamma_{2}$. In same time, we need to
estimate in a good way the dissipation phenomena in order to improve the
logarithmic decay rate. To this end, we apply some simple tools usually used
in the propagation of singularities \cite{AG}-\cite{CV}-\cite{R} in order to
link the number of reflections and the $s$ variable of the bicharacteristic
flow. In particular, we will work with the operator $O_{p}\left(  D\right)
=i\partial_{s}+h\left(  \Delta-\partial_{t}^{2}\right)  $ for $h\in\left(
0,h_{o}\right]  $. Observe that the product of four, mono-dimensional,
solutions of the Schr\"{o}dinger equation $i\partial_{s}\pm h\partial^{2}$ can
create a solution of $O_{p}\left(  D\right)  a\left(  x,t,s\right)  =0$. The
dispersive property of the linear Schr\"{o}dinger equation $i\partial_{s}\pm
h\partial^{2}$\ will be exploited.

\bigskip

The next section describes an interpolation inequality. In section 3, we give
the proof of Theorem. We close this paper with two Appendixes devoted to prove
some technical results.

\bigskip

\section{Interpolation inequality}

\bigskip

The purpose of this section is to establish the following inequality.

\bigskip

Theorem 2 .- \textit{Let }$T>0$\textit{. There exist }$C>0$\textit{ and
}$\gamma>1$\textit{ such that for any }$h\in\left(  0,h_{o}\right]  $\textit{
and initial data }$\left(  u_{0},u_{1}\right)  \in H^{2}\left(  \Omega\right)
\cap H_{0}^{1}\left(  \Omega\right)  \times H_{0}^{1}\left(  \Omega\right)
$\textit{, the solution }$u$\textit{ of }(\ref{1.2})\textit{ satisfies}%
\[%
\begin{array}
[c]{ll}%
\int_{\omega_{o}}\int_{0}^{T}\left|  \partial_{t}u\left(  x,t\right)  \right|
^{2}dxdt & \leq C\left(  \frac{1}{h}\right)  ^{\gamma}\left(  \int_{\Upsilon
}\int_{\left|  t\right|  \leq\gamma\left(  \frac{1}{h}\right)  ^{\gamma}%
}\left|  \partial_{\nu}u\left(  x,t\right)  \right|  ^{2}dxdt\right)
^{1/2}\left\|  \left(  u_{0},u_{1}\right)  \right\|  _{H_{0}^{1}\left(
\Omega\right)  \times L^{2}\left(  \Omega\right)  }\\
& +C\sqrt{h}~\left\|  \left(  u_{0},u_{1}\right)  \right\|  _{H^{2}\cap
H_{0}^{1}\left(  \Omega\right)  \times H_{0}^{1}\left(  \Omega\right)
}\left\|  \left(  u_{0},u_{1}\right)  \right\|  _{H_{0}^{1}\left(
\Omega\right)  \times L_{0}^{2}\left(  \Omega\right)  }\text{ .}%
\end{array}
\]

\bigskip

The rest of this section is devoted to the proof of Theorem 2. We begin to
introduce some weight functions inspired by \cite{Z}, with the properties of
localization, propagation and dispersion. Next, we make appear the Fourier
variables and introduce some Fourier integral operators dependent on the
number of reflections between $\Gamma_{1}$ and $\Gamma_{2}$. Their properties
at the boundary are analyzed. The process of propagation is then applied. Some
parameters are adequately chosen and that will complete the proof.

\bigskip

\subsection{The weight functions}

\bigskip

The hypothesis saying that $\Upsilon\subset\mathbb{R}^{2}\left\backslash
\left(  \left(  -m_{1},m_{1}\right)  \times\left(  -m_{2},m_{2}\right)
\right)  \right.  \times\mathbb{R}$ implies that for any $x_{o}\in
\overline{\omega_{o}}$, $B\left(  x_{o},r_{o}/2\right)  \cap\partial
\Omega=\emptyset$, where $B\left(  x_{o},r\right)  $ denotes the ball of
center $x_{o}$ and radius $r$. We introduce $\chi_{x_{o}}=\chi\in
C_{0}^{\infty}\left(  B\left(  x_{o},r_{o}/2\right)  \right)  $ be such that
$0\leq\chi\leq1$ and $\chi=1$ on $B\left(  x_{o},r_{o}/4\right)  $.

\bigskip

\noindent From any point $x_{o}\in\overline{\omega_{o}}$, we will localize
around $x_{o}$ and eventually be able to propagate some local regularity. To
this end, let $h\in\left(  0,h_{o}\right]  $ and let us define
\[
a\left(  x,t,s\right)  =\left(  \frac{1}{\left(  is+1\right)  ^{3/2}%
}~e^{-\frac{1}{4h}\frac{x^{2}}{is+1}}\right)  \left(  \frac{1}{\sqrt{-ihs+1}%
}~e^{-\frac{1}{4}\frac{t^{2}}{-ihs+1}}\right)  \text{ ,}%
\]%

\[
a_{o}\left(  x,t\right)  =a\left(  x-x_{o},t,0\right)  \text{ and }%
\varphi\left(  x,t\right)  =\chi\left(  x\right)  a_{o}\left(  x,t\right)
\text{ .}%
\]

\noindent We get the following identities%

\[
\left(  i\partial_{s}+h\left(  \Delta-\partial_{t}^{2}\right)  \right)
a\left(  x,t,s\right)  =0\text{\quad}\forall\left(  x,t,s\right)  \in
\Omega\times\mathbb{R}\times\left(  0,+\infty\right)  \text{ ,}%
\]%

\begin{equation}
\left|  a\left(  x,t,s\right)  \right|  =\frac{1}{\left(  \sqrt{s^{2}%
+1}\right)  ^{3/2}}\frac{1}{\left(  \sqrt{\left(  hs\right)  ^{2}+1}\right)
^{1/2}}~e^{-\frac{x^{2}}{4h}\frac{1}{s^{2}+1}}~e^{-\frac{t^{2}}{4}\frac
{1}{\left(  hs\right)  ^{2}+1}}\text{ .} \tag{2.1}\label{2.1}%
\end{equation}

\noindent Now we use such weight functions $a_{o}$ and $\varphi$ with $u$ the
solution of the wave equation (\ref{1.2}) as follows. By integrations by
parts,
\[%
\begin{array}
[c]{ll}%
& \quad\int_{\Omega\times\mathbb{R}}\chi\left(  x\right)  \left|
a_{o}\partial_{t}u\left(  x,t\right)  \right|  ^{2}dxdt\\
& =-\int_{\Omega\times\mathbb{R}}\chi\left(  x\right)  \partial_{t}\left(
\left|  a_{o}\left(  x,t\right)  \right|  ^{2}\right)  \frac{1}{2}\partial
_{t}\left(  \left|  u\left(  x,t\right)  \right|  ^{2}\right)  dxdt-\int
_{\Omega\times\mathbb{R}}\chi\left(  x\right)  \left|  a_{o}\left(
x,t\right)  \right|  ^{2}\partial_{t}^{2}u\left(  x,t\right)  u\left(
x,t\right)  dxdt\\
& =\int_{\Omega\times\mathbb{R}}\chi\left(  x\right)  \frac{1}{2}\partial
_{t}^{2}\left(  \left|  a_{o}\left(  x,t\right)  \right|  ^{2}\right)  \left|
u\left(  x,t\right)  \right|  ^{2}dxdt-\int_{\Omega\times\mathbb{R}}%
\chi\left(  x\right)  \left|  a_{o}\left(  x,t\right)  \right|  ^{2}%
\partial_{t}^{2}u\left(  x,t\right)  u\left(  x,t\right)  dxdt\text{ .}%
\end{array}
\]

\noindent As $\partial_{t}^{2}\left(  \left|  a_{o}\left(  x,t\right)
\right|  ^{2}\right)  =-\left|  a_{o}\left(  x,t\right)  \right|  ^{2}%
+t^{2}\left|  a_{o}\left(  x,t\right)  \right|  ^{2}$ and $t^{2}\left|
a_{o}\left(  x,t\right)  \right|  ^{2}\leq4\left|  a_{o}\left(  x,t/\sqrt
{2}\right)  \right|  ^{2}$, we have
\begin{equation}%
\begin{array}
[c]{ll}%
& \quad\int_{\Omega\times\mathbb{R}}\chi\left(  x\right)  \left|
a_{o}\partial_{t}u\left(  x,t\right)  \right|  ^{2}dxdt\\
& \leq2\int_{\Omega\times\mathbb{R}}a_{o}\left(  x,t/\sqrt{2}\right)
\varphi\left(  x,t/\sqrt{2}\right)  \left|  u\left(  x,t\right)  \right|
^{2}dxdt+\left|  \int_{\Omega\times\mathbb{R}}a_{o}\left(  x,t\right)
\varphi\left(  x,t\right)  \partial_{t}^{2}u\left(  x,t\right)  u\left(
x,t\right)  dxdt\right|  \text{ .}%
\end{array}
\tag{2.2}\label{2.2}%
\end{equation}

\bigskip

\noindent Let $\left\{  G_{i};i\in I\right\}  $ be a family of open sets
covering $\omega_{o}$, i.e., $\overline{\omega_{o}}\subset\bigcup\limits_{i\in
I}G_{i}$, such that $G_{i}=\left\{  \left|  x-x_{o}^{i}\right|  \leq2\sqrt
{h}\right\}  $ where $\left\{  x_{o}^{i}\right\}  _{i\in I}\in\overline
{\omega_{o}}$ and $I$ is a countable set such that the number of elements of
$I$ is bounded by $\frac{c_{o}}{h\sqrt{h}}$ for some constant $c_{o}>0$
independent of $h\in\left(  0,h_{o}\right]  $. Consequently,
\begin{equation}%
\begin{array}
[c]{ll}%
\int_{\omega_{o}\times\left(  0,T\right)  }\left|  \partial_{t}u\left(
x,t\right)  \right|  ^{2}dxdt & \leq e^{\frac{1}{2}T^{2}}\int_{\omega
_{o}\times\left(  0,T\right)  }e^{-\frac{1}{2}t^{2}}\left|  \partial
_{t}u\left(  x,t\right)  \right|  ^{2}dxdt\\
& \leq e^{\frac{1}{2}T^{2}+2}\sum\limits_{i\in I}\int_{G_{i}\times\mathbb{R}%
}\chi_{x_{o}^{i}}\left(  x\right)  \left|  a\left(  x-x_{o}^{i},t,0\right)
\partial_{t}u\left(  x,t\right)  \right|  ^{2}dxdt\\
& \leq e^{\frac{1}{2}T^{2}+2}\sum\limits_{i\in I}\int_{\Omega\times\mathbb{R}%
}\chi_{x_{o}^{i}}\left(  x\right)  \left|  a\left(  x-x_{o}^{i},t,0\right)
\partial_{t}u\left(  x,t\right)  \right|  ^{2}dxdt
\end{array}
\tag{2.3}\label{2.3}%
\end{equation}

\noindent and we will search to bound $\int_{\Omega\times\mathbb{R}}%
\chi_{x_{o}^{i}}\left(  x\right)  \left|  a\left(  x-x_{o}^{i},t,0\right)
\partial_{t}u\left(  x,t\right)  \right|  ^{2}dxdt$ by a suitable term
$E_{h}\left(  u\right)  $ independent of $i$ in order to get
\[
\int_{\omega_{o}\times\left(  0,T\right)  }\left|  \partial_{t}u\left(
x,t\right)  \right|  ^{2}dxdt\leq\frac{C_{o}E_{h}\left(  u\right)  }{h\sqrt
{h}}%
\]

\noindent for some $C_{o}>0$ independent of $u$ and $h$. To this end, we will
first study in the next subsections the second term of the second member of
(\ref{2.2}),%
\[
\int_{\Omega\times\mathbb{R}}a_{o}\left(  x,t\right)  \varphi\left(
x,t\right)  f\left(  x,t\right)  u\left(  x,t\right)  dxdt\text{ when
}f=\partial_{t}^{2}u\text{.}%
\]

\bigskip

\subsection{The Fourier variables}

\bigskip

\noindent Denote%
\[
\widehat{\varphi f}\left(  \xi,\tau\right)  =\int_{\Omega\times\mathbb{R}%
}e^{-i\left(  x\xi+t\tau\right)  }~\varphi\left(  x,t\right)  f\left(
x,t\right)  ~dxdt\text{ ,}%
\]

\noindent then for any $\left(  x,t\right)  \in\Omega\times\mathbb{R}$,
\[
a_{o}\left(  x,t\right)  \varphi\left(  x,t\right)  f\left(  x,t\right)
=a_{o}\left(  x,t\right)  \frac{1}{\left(  2\pi\right)  ^{4}}\int
_{\mathbb{R}^{4}}e^{i\left(  x\xi+t\tau\right)  }~\widehat{\varphi f}\left(
\xi,\tau\right)  ~d\xi d\tau\text{ .}%
\]

\noindent Let $\lambda\geq1$. We cut the integral over $\tau\in\mathbb{R}$
into two parts, $\left\{  \left|  \tau\right|  \geq\lambda\right\}  $ and
$\left\{  \left|  \tau\right|  <\lambda\right\}  $. Next, for $\left\{
\left|  \tau\right|  <\lambda\right\}  $ and $\xi=\left(  \xi_{1},\xi_{2}%
,\xi_{3}\right)  $, the integral over $\xi_{3}\in\mathbb{R}$ is divided as
follows. Denote $\left(  2\mathbb{Z}+1\right)  =\left\{  2n+1\left\backslash
n\in\mathbb{Z}\right.  \right\}  $, then%
\[%
\begin{array}
[c]{ll}%
a_{o}\left(  x,t\right)  \varphi\left(  x,t\right)  f\left(  x,t\right)  &
=a_{o}\left(  x,t\right)  \frac{1}{\left(  2\pi\right)  ^{4}}\sum
\limits_{\xi_{o3}\in\left(  2\mathbb{Z}+1\right)  }\int_{\mathbb{R}^{2}}%
\int_{\xi_{o3}-1}^{\xi_{o3}+1}\int_{\left|  \tau\right|  <\lambda}e^{i\left(
x\xi+t\tau\right)  }~\widehat{\varphi f}\left(  \xi,\tau\right)  ~d\xi d\tau\\
& \quad+a_{o}\left(  x,t\right)  \frac{1}{\left(  2\pi\right)  ^{4}}%
\int_{\mathbb{R}^{3}}\int_{\left|  \tau\right|  \geq\lambda}e^{i\left(
x\xi+t\tau\right)  }~\widehat{\varphi f}\left(  \xi,\tau\right)  ~d\xi
d\tau\text{ .}%
\end{array}
\]

\noindent Consequently,
\begin{equation}%
\begin{array}
[c]{ll}%
& \quad\int_{\Omega\times\mathbb{R}}a_{o}\left(  x,t\right)  \varphi\left(
x,t\right)  f\left(  x,t\right)  u\left(  x,t\right)  dxdt-R_{0}\\
& =\int_{\Omega\times\mathbb{R}}a_{o}\left(  x,t\right)  \frac{1}{\left(
2\pi\right)  ^{4}}\sum\limits_{\xi_{o3}\in\left(  2\mathbb{Z}+1\right)  }%
\int_{\mathbb{R}^{2}}\int_{\xi_{o3}-1}^{\xi_{o3}+1}\int_{\left|  \tau\right|
<\lambda}e^{i\left(  x\xi+t\tau\right)  }~\widehat{\varphi f}\left(  \xi
,\tau\right)  ~d\xi d\tau~u\left(  x,t\right)  dxdt
\end{array}
\tag{2.4}\label{2.4}%
\end{equation}

\noindent where%
\[
R_{0}=\int_{\Omega\times\mathbb{R}}a_{o}\left(  x,t\right)  \frac{1}{\left(
2\pi\right)  ^{4}}\int_{\mathbb{R}^{3}}\int_{\left|  \tau\right|  \geq\lambda
}e^{i\left(  x\xi+t\tau\right)  }~\widehat{\varphi f}\left(  \xi,\tau\right)
~d\xi d\tau~u\left(  x,t\right)  ~dxdt\text{ .}%
\]

\noindent Moreover, it holds%
\begin{equation}
R_{0}\leq c~\sqrt{\frac{1}{\lambda}}\sqrt{\mathcal{G}\left(  u,0\right)
}\sqrt{\mathcal{G}\left(  \partial_{t}u,0\right)  }\text{ .} \tag{2.5}%
\label{2.5}%
\end{equation}

\noindent The proof of (\ref{2.5}) is given in Appendix A.

\bigskip

\subsection{The Fourier integral operators}

\bigskip

\noindent Denote $x=\left(  x_{1},x_{2},x_{3}\right)  \in\Omega$, $\xi=\left(
\xi_{1},\xi_{2},\xi_{3}\right)  $ with $\left(  \xi_{1},\xi_{2}\right)
\in\mathbb{R}^{2}$. Let $\left(  x_{o},\xi_{o3}\right)  =\left(  x_{o1}%
,x_{o2},x_{o3},\xi_{o3}\right)  \in\overline{\omega_{o}}\times\left(
2\mathbb{Z}+1\right)  $, we introduce for all $s\geq0$ and $n\in\mathbb{Z}$
\[%
\begin{array}
[c]{ll}%
& \quad A_{s}^{n}\left(  x_{o},\xi_{o3}\right)  f\left(  x,t\right) \\
& =\frac{\left(  -1\right)  ^{n}}{\left(  2\pi\right)  ^{4}}\int
_{\mathbb{R}^{2}}\int_{\xi_{o3}-1}^{\xi_{o3}+1}\int_{\left|  \tau\right|
<\lambda}e^{i\left(  x_{1}\xi_{1}+x_{2}\xi_{2}+t\tau\right)  }~e^{i\left[
\left(  -1\right)  ^{n}x_{3}+2n\frac{\xi_{o3}}{\left|  \xi_{o3}\right|  }%
\rho\right]  \xi_{3}}~e^{-i\left(  \xi^{2}-\tau^{2}\right)  hs}~\widehat
{\varphi f}\left(  \xi,\tau\right) \\
& \qquad a\left(  x_{1}-x_{o1}-2\xi_{1}hs,x_{2}-x_{o2}-2\xi_{2}hs,x_{3}%
-\left(  -1\right)  ^{n}\left[  -2n\frac{\xi_{o3}}{\left|  \xi_{o3}\right|
}\rho+x_{o3}+2\xi_{3}hs\right]  ,t+2\tau hs,s\right)  d\xi d\tau\text{ .}%
\end{array}
\]

\bigskip

\noindent Let $\left(  P,Q\right)  \in\mathbb{N}^{2}$. Consider the solution
\[
A_{s,P,Q}\left(  x_{o},\xi_{o3}\right)  f\left(  x,t\right)  =\sum
\limits_{n=-2Q}^{2P+1}A_{s}^{n}\left(  x_{o},\xi_{o3}\right)  f\left(
x,t\right)  \text{ .}%
\]

\noindent One can check that for any $\left(  x_{o},\xi_{o3},P,Q\right)
\in\overline{\omega_{o}}\times\left(  2\mathbb{Z}+1\right)  \times
\mathbb{N}^{2}$,
\begin{equation}
\left(  i\partial_{s}+h\left(  \Delta-\partial_{t}^{2}\right)  \right)
A_{s,P,Q}\left(  x_{o},\xi_{o3}\right)  f\left(  x,t\right)  =0\text{\quad
}\forall\left(  x,t,s\right)  \in\Omega\times\mathbb{R}\times\left(
0,+\infty\right)  \text{ .} \tag{2.6}\label{2.6}%
\end{equation}

\bigskip

\subsection{At $s=0$}

\bigskip

Since%
\[
A_{0}^{0}\left(  x_{o},\xi_{o3}\right)  f\left(  x,t\right)  =a_{o}\left(
x,t\right)  \frac{1}{\left(  2\pi\right)  ^{4}}\int_{\mathbb{R}^{2}}\int
_{\xi_{o3}-1}^{\xi_{o3}+1}\int_{\left|  \tau\right|  <\lambda}e^{i\left(
x\xi+t\tau\right)  }~\widehat{\varphi f}\left(  \xi,\tau\right)  ~d\xi
d\tau\text{ ,}%
\]

\noindent we then obtain from (\ref{2.4}) that
\[
\int_{\Omega\times\mathbb{R}}a_{o}\left(  x,t\right)  \varphi\left(
x,t\right)  f\left(  x,t\right)  u\left(  x,t\right)  dxdt-R_{0}=\int
_{\Omega\times\mathbb{R}}\sum\limits_{\xi_{o3}\in\left(  2\mathbb{Z}+1\right)
}A_{0}^{0}\left(  x_{o},\xi_{o3}\right)  f\left(  x,t\right)  ~u\left(
x,t\right)  dxdt\text{ .}%
\]

\noindent Observe that
\[
A_{0,P,Q}\left(  x_{o},\xi_{o3}\right)  f\left(  x,t\right)  =A_{0}^{0}\left(
x_{o},\xi_{o3}\right)  f\left(  x,t\right)  +\sum\limits_{n=1}^{2P+1}A_{0}%
^{n}\left(  x_{o},\xi_{o3}\right)  f\left(  x,t\right)  +\sum\limits_{-2Q\leq
n\leq-1}A_{0}^{n}\left(  x_{o},\xi_{o3}\right)  f\left(  x,t\right)  \text{ ,}%
\]

\noindent with the convention that for $Q=0$, $\sum\limits_{-2Q\leq n\leq
-1}A_{0}^{n}\left(  x_{o},\xi_{o3}\right)  f\left(  x,t\right)  =0$. So, we
deduce that%
\begin{equation}
\int_{\Omega\times\mathbb{R}}a_{o}\left(  x,t\right)  \varphi\left(
x,t\right)  f\left(  x,t\right)  u\left(  x,t\right)  dxdt-R_{0}-R_{1}%
=\int_{\Omega\times\mathbb{R}}\sum\limits_{\xi_{o3}\in\left(  2\mathbb{Z}%
+1\right)  }A_{0,P,Q}\left(  x_{o},\xi_{o3}\right)  f\left(  x,t\right)
~u\left(  x,t\right)  dxdt \tag{2.7}\label{2.7}%
\end{equation}

\noindent where%
\[
R_{1}=-\int_{\Omega\times\mathbb{R}}\sum\limits_{\xi_{o3}\in\left(
2\mathbb{Z}+1\right)  }\left[  \sum\limits_{n=1}^{2P+1}A_{0}^{n}\left(
x_{o},\xi_{o3}\right)  f\left(  x,t\right)  +\sum\limits_{-2Q\leq n\leq
-1}A_{0}^{n}\left(  x_{o},\xi_{o3}\right)  f\left(  x,t\right)  \right]
~u\left(  x,t\right)  dxdt\text{ .}%
\]

\noindent We estimate $R_{1}$ uniformly with respect to $\left(  P,Q\right)  $
as follows.%
\[%
\begin{array}
[c]{ll}%
R_{1} & \leq\int_{\Omega\times\mathbb{R}}\sum\limits_{\xi_{o3}\in\left(
2\mathbb{Z}+1\right)  }\sum\limits_{n\in\mathbb{Z}\left\backslash \left\{
0\right\}  \right.  }\left|  A_{0}^{n}\left(  x_{o},\xi_{o3}\right)  f\left(
x,t\right)  \right|  ~\left|  u\left(  x,t\right)  \right|  dxdt\\
& \leq\int_{\Omega\times\mathbb{R}}\sum\limits_{\xi_{o3}\in\left(
2\mathbb{Z}+1\right)  }\frac{1}{\left(  2\pi\right)  ^{4}}\int_{\mathbb{R}%
^{2}}\int_{\xi_{o3}-1}^{\xi_{o3}+1}\int_{\left|  \tau\right|  <\lambda}\left|
\widehat{\varphi f}\left(  \xi,\tau\right)  \right|  ~d\xi d\tau\\
& \quad\quad\quad\sum\limits_{n\in\mathbb{Z}\left\backslash \left\{
0\right\}  \right.  }a\left(  x_{1}-x_{o1},x_{2}-x_{o2},\left(  -1\right)
^{n}x_{3}+2n\frac{\xi_{o3}}{\left|  \xi_{o3}\right|  }\rho-x_{o3},t,0\right)
~\left|  u\left(  x,t\right)  \right|  dxdt\\
& \leq\frac{1}{\left(  2\pi\right)  ^{4}}\int_{\mathbb{R}^{3}}\int_{\left|
\tau\right|  <\lambda}\left|  \widehat{\varphi f}\left(  \xi,\tau\right)
\right|  ~d\xi d\tau\\
& \quad\quad\quad\int_{\Omega\times\mathbb{R}}\sum\limits_{n\in\mathbb{Z}%
\left\backslash \left\{  0\right\}  \right.  }a\left(  x_{1}-x_{o1}%
,x_{2}-x_{o2},\left(  -1\right)  ^{n}x_{3}+2n\rho-x_{o3},t,0\right)  ~\left|
u\left(  x,t\right)  \right|  dxdt\text{ .}%
\end{array}
\]

\noindent Now, notice that
\[%
\begin{array}
[c]{cc}%
& \quad a\left(  x_{1}-x_{o1},x_{2}-x_{o2},\left(  -1\right)  ^{n}%
x_{3}+2n\frac{\xi_{o3}}{\left|  \xi_{o3}\right|  }\rho-x_{o3},t,0\right) \\
& =e^{-\frac{1}{4h}\left[  \left(  x_{1}-x_{o1}\right)  ^{2}+\left(
x_{2}-x_{o2}\right)  ^{2}+\left(  \left(  -1\right)  ^{n}x_{3}+2n\frac
{\xi_{o3}}{\left|  \xi_{o3}\right|  }\rho-x_{o3}\right)  ^{2}\right]
}~e^{-\frac{t^{2}}{4}}\text{ .}%
\end{array}
\]

\noindent The hypothesis saying that $\Upsilon\subset\mathbb{R}^{2}%
\left\backslash \left(  \left(  -m_{1},m_{1}\right)  \times\left(
-m_{2},m_{2}\right)  \right)  \right.  \times\mathbb{R}$ implies that for any
$x_{o}\in\overline{\omega_{o}}$ and $x=\left(  x_{1},x_{2},x_{3}\right)
\in\Omega$ such that $\left(  x_{1},x_{2}\right)  \notin\left[  -m_{1}%
+r_{o}/2,m_{1}-r_{o}/2\right]  \times\left[  -m_{2}+r_{o}/2,m_{2}%
-r_{o}/2\right]  $, we have $\left(  x_{1}-x_{o1}\right)  ^{2}+\left(
x_{2}-x_{o2}\right)  ^{2}\geq\left(  r_{o}/2\right)  ^{2}$, but such
hypothesis also implies that for any $x_{o}\in\overline{\omega_{o}}$ and
$x=\left(  x_{1},x_{2},x_{3}\right)  \in\Omega$ such that $\left(  x_{1}%
,x_{2}\right)  \in\left[  -m_{1}+r_{o}/2,m_{1}-r_{o}/2\right]  \times\left[
-m_{2}+r_{o}/2,m_{2}-r_{o}/2\right]  $, we get $x_{3}\in\left[  -\rho
,\rho\right]  $ and therefore $\left|  \left(  -1\right)  ^{n}x_{3}%
+2n\frac{\xi_{o3}}{\left|  \xi_{o3}\right|  }\rho-x_{o3}\right|  \geq\frac
{3}{4}\rho$ for any $n\in\mathbb{Z}\left\backslash \left\{  0\right\}
\right.  $. So, for any $n\in\mathbb{Z}\left\backslash \left\{  0\right\}
\right.  $%
\[
a\left(  x_{1}-x_{o1},x_{2}-x_{o2},\left(  -1\right)  ^{n}x_{3}+2n\frac
{\xi_{o3}}{\left|  \xi_{o3}\right|  }\rho-x_{o3},t,0\right)  \leq e^{-\frac
{c}{h}}e^{-\frac{1}{8h}\left(  \left(  -1\right)  ^{n}x_{3}+2n\frac{\xi_{o3}%
}{\left|  \xi_{o3}\right|  }\rho-x_{o3}\right)  ^{2}}e^{-\frac{t^{2}}{4}%
}\text{ .}%
\]

\noindent It follows that%
\[
\int_{\Omega\times\mathbb{R}}\sum\limits_{n\in\mathbb{Z}\left\backslash
\left\{  0\right\}  \right.  }a\left(  x_{1}-x_{o1},x_{2}-x_{o2},\left(
-1\right)  ^{n}x_{3}+2n\rho-x_{o3},t,0\right)  ~\left|  u\left(  x,t\right)
\right|  dxdt\leq c~e^{-\frac{c}{h}}\sqrt{\mathcal{G}\left(  u,0\right)  }%
\]

\noindent Now it remains to compute $\int_{\mathbb{R}^{3}}\int_{\left|
\tau\right|  <\lambda}\left|  \widehat{\varphi f}\left(  \xi,\tau\right)
\right|  ~d\xi d\tau$. We can check that there exists $\gamma>1$ such that
\begin{equation}
\int_{\mathbb{R}^{3}}\int_{\left|  \tau\right|  <\lambda}\left|
\widehat{\varphi f}\left(  \xi,\tau\right)  \right|  ~d\xi d\tau\leq c\left(
\frac{\lambda}{h}\right)  ^{\gamma}\sqrt{\mathcal{G}\left(  u,0\right)
}\text{ .} \tag{2.8}\label{2.8}%
\end{equation}

\noindent The proof of (\ref{2.8}) is given in Appendix A.

\bigskip

\noindent We conclude that
\begin{equation}
R_{1}\leq c\left(  \frac{\lambda}{h}\right)  ^{\gamma}e^{-\frac{c}{h}%
}\mathcal{G}\left(  u,0\right)  \text{ .} \tag{2.9}\label{2.9}%
\end{equation}

\bigskip

\subsection{On the boundary $\left\{  x_{3}=\pm\rho\right\}  $}

\bigskip

Since $a\left(  x_{1},x_{2},x_{3},t,s\right)  =a\left(  x_{1},x_{2}%
,-x_{3},t,s\right)  $, we get the following identity%
\begin{equation}
A_{s}^{n}\left(  x_{o},\xi_{o3}\right)  f\left(  x_{1},x_{2},\left(
-1\right)  ^{n}\frac{\xi_{o3}}{\left|  \xi_{o3}\right|  }\rho,t\right)
=-A_{s}^{n+1}\left(  x_{o},\xi_{o3}\right)  f\left(  x_{1},x_{2},\left(
-1\right)  ^{n}\frac{\xi_{o3}}{\left|  \xi_{o3}\right|  }\rho,t\right)
\quad\forall n\in\mathbb{Z}\text{ .} \tag{2.10}\label{2.10}%
\end{equation}

\noindent Thus,
\[
A_{s}^{2n}\left(  x_{o},\xi_{o3}\right)  f\left(  x_{1},x_{2},\frac{\xi_{o3}%
}{\left|  \xi_{o3}\right|  }\rho,t\right)  =-A_{s}^{2n+1}\left(  x_{o}%
,\xi_{o3}\right)  f\left(  x_{1},x_{2},\frac{\xi_{o3}}{\left|  \xi
_{o3}\right|  }\rho,t\right)  \quad\forall n\in\mathbb{Z}\text{ ,}%
\]

\noindent so that
\[%
\begin{array}
[c]{ll}%
A_{s,P,Q}\left(  x_{o},\xi_{o3}\right)  f\left(  x_{1},x_{2},\frac{\xi_{o3}%
}{\left|  \xi_{o3}\right|  }\rho,t\right)  & =\sum\limits_{n=-Q}^{P}\left[
A_{s}^{2n}\left(  x_{o},\xi_{o3}\right)  f\left(  x,t\right)  +A_{s}%
^{2n+1}\left(  x_{o},\xi_{o3}\right)  f\left(  x,t\right)  \right] \\
& =0\quad\forall s\geq0\text{ .}%
\end{array}
\]

\noindent Also, (\ref{2.10}) implies%
\[
A_{s}^{2n+1}\left(  x_{o},\xi_{o3}\right)  f\left(  x_{1},x_{2},-\frac
{\xi_{o3}}{\left|  \xi_{o3}\right|  }\rho,t\right)  =-A_{s}^{2n+2}\left(
x_{o},\xi_{o3}\right)  f\left(  x_{1},x_{2},-\frac{\xi_{o3}}{\left|  \xi
_{o3}\right|  }\rho,t\right)  \quad\forall n\in\mathbb{Z}\text{ ,}%
\]

\noindent therefore%
\[%
\begin{array}
[c]{ll}%
& \quad A_{s,P,Q}\left(  x_{o},\xi_{o3}\right)  f\left(  x_{1},x_{2}%
,-\frac{\xi_{o3}}{\left|  \xi_{o3}\right|  }\rho,t\right) \\
& =A_{s}^{-2Q}\left(  x_{o},\xi_{o3}\right)  f\left(  x_{1},x_{2},-\frac
{\xi_{o3}}{\left|  \xi_{o3}\right|  }\rho,t\right)  +A_{s}^{2P+1}\left(
x_{o},\xi_{o3}\right)  f\left(  x_{1},x_{2},-\frac{\xi_{o3}}{\left|  \xi
_{o3}\right|  }\rho,t\right)  \quad\forall s\geq0\text{ .}%
\end{array}
\]

\bigskip

\subsection{The key identity}

\bigskip

\noindent By multiplying (\ref{2.6}) by $u\left(  x,t\right)  $ and
integrating by parts over $\Omega\times\mathbb{R}\times\left[  0,L\right]  $,
we have that for all $\left(  x_{o},\xi_{o3}\right)  \in\overline{\omega_{o}%
}\times\left(  2\mathbb{Z}+1\right)  $, for all $\left(  P,Q\right)
\in\mathbb{N}^{2}$ and all $L>0$,%
\[%
\begin{array}
[c]{ll}%
\int_{\Omega\times\mathbb{R}}A_{0,P,Q}\left(  x_{o},\xi_{o3}\right)  f\left(
x,t\right)  u\left(  x,t\right)  dxdt & =\int_{\Omega\times\mathbb{R}%
}A_{L,P,Q}\left(  x_{o},\xi_{o3}\right)  f\left(  x,t\right)  u\left(
x,t\right)  dxdt\\
& \quad+ih\int_{0}^{L}\int_{\partial\Omega\times\mathbb{R}}A_{s,P,Q}\left(
x_{o},\xi_{o3}\right)  f\left(  x,t\right)  \partial_{\nu}u\left(  x,t\right)
dxdtds\text{ .}%
\end{array}
\]

\noindent Consequently, combining the later equality with (\ref{2.7}), we have
the following key identity%
\begin{equation}%
\begin{array}
[c]{ll}%
& \quad\int_{\Omega\times\mathbb{R}}a_{o}\left(  x,t\right)  \varphi\left(
x,t\right)  f\left(  x,t\right)  u\left(  x,t\right)  dxdt-R_{0}-R_{1}\\
& =\int_{\Omega\times\mathbb{R}}\sum\limits_{\xi_{o3}\in\left(  2\mathbb{Z}%
+1\right)  }A_{L,P\left(  \xi_{o3}\right)  ,Q\left(  \xi_{o3}\right)  }\left(
x_{o},\xi_{o3}\right)  f\left(  x,t\right)  u\left(  x,t\right)  dxdt\\
& \quad+ih\int_{0}^{L}\int_{\partial\Omega\times\mathbb{R}}\sum\limits_{\xi
_{o3}\in\left(  2\mathbb{Z}+1\right)  }A_{s,P\left(  \xi_{o3}\right)
,Q\left(  \xi_{o3}\right)  }\left(  x_{o},\xi_{o3}\right)  f\left(
x,t\right)  \partial_{\nu}u\left(  x,t\right)  dxdtds
\end{array}
\tag{2.11}\label{2.11}%
\end{equation}

\noindent for any $L>0$ and $\left(  P\left(  \xi_{o3}\right)  ,Q\left(
\xi_{o3}\right)  \right)  \in\mathbb{N}^{2}$. $L$ will be taken large enough
in order that the first term in the second member of (\ref{2.11}) is
polynomially small like $\frac{c}{\sqrt{L}}$ and this uniformly with respect
to $\left(  \xi_{o3},P\left(  \xi_{o3}\right)  ,Q\left(  \xi_{o3}\right)
\right)  $ by using the dispersion of (\ref{2.1}). Next, for each $\xi_{o3}%
\in\left(  2\mathbb{Z}+1\right)  $, $\left(  P\left(  \xi_{o3}\right)
,Q\left(  \xi_{o3}\right)  \right)  $ will be chosen dependent on $\left(
\xi_{o3},L\right)  $ in order that the second term in the second member of
(\ref{2.11}) is exponentially small with respect to $h$ on the boundary
$\Gamma_{1}\cup\Gamma_{2}$.

\bigskip

\subsection{The internal term}

\bigskip

In this subsection, we study the internal term appearing in (\ref{2.11})%
\[
\int_{\Omega\times\mathbb{R}}\sum\limits_{\xi_{o3}\in\left(  2\mathbb{Z}%
+1\right)  }A_{L,P,Q}\left(  x_{o},\xi_{o3}\right)  f\left(  x,t\right)
u\left(  x,t\right)  dxdt\text{ .}%
\]

\noindent First, we have a uniform bound with respect to $\left(  P,Q\right)
$%
\begin{equation}%
\begin{array}
[c]{ll}%
& \quad\int_{\Omega\times\mathbb{R}}\sum\limits_{\xi_{o3}\in\left(
2\mathbb{Z}+1\right)  }A_{L,P,Q}\left(  x_{o},\xi_{o3}\right)  f\left(
x,t\right)  u\left(  x,t\right)  dxdt\\
& =\int_{\Omega\times\mathbb{R}}\sum\limits_{\xi_{o3}\in\left(  2\mathbb{Z}%
+1\right)  }\left[  \sum\limits_{n=-2Q}^{2P+1}A_{L}^{n}\left(  x_{o},\xi
_{o3}\right)  f\left(  x,t\right)  \right]  u\left(  x,t\right)  dxdt\\
& \leq\sum\limits_{\xi_{o3}\in\left(  2\mathbb{Z}+1\right)  }\sum
\limits_{n\in\mathbb{Z}}\left|  \int_{\Omega\times\mathbb{R}}A_{L}^{n}\left(
x_{o},\xi_{o3}\right)  f\left(  x,t\right)  u\left(  x,t\right)  dxdt\right|
\text{ .}%
\end{array}
\tag{2.12}\label{2.12}%
\end{equation}

\noindent Recall that
\[%
\begin{array}
[c]{ll}%
& \quad A_{L}^{n}\left(  x_{o},\xi_{o3}\right)  f\left(  x,t\right) \\
& =\frac{1}{\left(  iL+1\right)  ^{3/2}}\frac{\left(  -1\right)  ^{n}}{\left(
2\pi\right)  ^{4}}\int_{\mathbb{R}^{2}}\int_{\xi_{o3}-1}^{\xi_{o3}+1}%
\int_{\left|  \tau\right|  <\lambda}e^{i\left(  x_{1}\xi_{1}+x_{2}\xi
_{2}+t\tau\right)  }~e^{i\left[  \left(  -1\right)  ^{n}x_{3}+2n\frac{\xi
_{o3}}{\left|  \xi_{o3}\right|  }\rho\right]  \xi_{3}}~e^{-i\left(  \xi
^{2}-\tau^{2}\right)  hL}~\widehat{\varphi f}\left(  \xi,\tau\right) \\
& \qquad\quad e^{-\frac{1}{4h}\frac{\left(  x_{1}-x_{o1}-2\xi_{1}hL\right)
^{2}}{iL+1}}~e^{-\frac{1}{4h}\frac{\left(  x_{2}-x_{o2}-2\xi_{2}hL\right)
^{2}}{iL+1}}~e^{-\frac{1}{4h}\frac{\left(  \left(  -1\right)  ^{n}%
x_{3}+2n\frac{\xi_{o3}}{\left|  \xi_{o3}\right|  }\rho-x_{o3}-2\xi
_{3}hL\right)  ^{2}}{iL+1}}\\
& \qquad\quad\left(  \frac{1}{\sqrt{-ihL+1}}~e^{-\frac{1}{4}\frac{\left(
t+2\tau hL\right)  ^{2}}{-ihL+1}}\right)  ~d\xi d\tau\text{ }%
\end{array}
\]%

\[%
\begin{array}
[c]{ll}%
u\left(  x,t\right)  & =\sum\limits_{j\geq1}\left[  \frac{b_{j}^{0}}{2}\left(
e^{it\sqrt{\mu_{j}}}+e^{-it\sqrt{\mu_{j}}}\right)  +\frac{b_{j}^{1}}%
{2i\sqrt{\mu_{j}}}\left(  e^{it\sqrt{\mu_{j}}}-e^{-it\sqrt{\mu_{j}}}\right)
\right]  ~\ell_{j}\left(  x\right) \\
& =\sum\limits_{j\geq1}\left(  \frac{b_{j}^{0}}{2}+\frac{b_{j}^{1}}%
{2i\sqrt{\mu_{j}}}\right)  e^{it\sqrt{\mu_{j}}}~\ell_{j}\left(  x\right)
+\sum\limits_{j\geq1}\left(  \frac{b_{j}^{0}}{2}-\frac{b_{j}^{1}}{2i\sqrt
{\mu_{j}}}\right)  e^{-it\sqrt{\mu_{j}}}~\ell_{j}\left(  x\right)  \text{ .}%
\end{array}
\]

\noindent Therefore, we write%
\begin{equation}%
\begin{array}
[c]{ll}%
& \quad\int_{\Omega\times\mathbb{R}}A_{L}^{n}\left(  x_{o},\xi_{o3}\right)
f\left(  x,t\right)  u\left(  x,t\right)  dxdt\\
& =\frac{1}{\left(  iL+1\right)  ^{3/2}}\frac{\left(  -1\right)  ^{n}}{\left(
2\pi\right)  ^{4}}\sum\limits_{j\geq1}\int_{\Omega}\ell_{j}\left(  x\right) \\
& \qquad\int_{\mathbb{R}^{2}}\int_{\xi_{o3}-1}^{\xi_{o3}+1}\int_{\left|
\tau\right|  <\lambda}e^{i\left(  x_{1}\xi_{1}+x_{2}\xi_{2}\right)
}~e^{i\left[  \left(  -1\right)  ^{n}x_{3}+2n\frac{\xi_{o3}}{\left|  \xi
_{o3}\right|  }\rho\right]  \xi_{3}}~e^{-i\left(  \xi^{2}-\tau^{2}\right)
hL}~\widehat{\varphi f}\left(  \xi,\tau\right)  ~d\xi\\
& \qquad\quad e^{-\frac{1}{4h}\frac{\left(  x_{1}-x_{o1}-2\xi_{1}hL\right)
^{2}}{iL+1}}~e^{-\frac{1}{4h}\frac{\left(  x_{2}-x_{o2}-2\xi_{2}hL\right)
^{2}}{iL+1}}~e^{-\frac{1}{4h}\frac{\left(  \left(  -1\right)  ^{n}%
x_{3}+2n\frac{\xi_{o3}}{\left|  \xi_{o3}\right|  }\rho-x_{o3}-2\xi
_{3}hL\right)  ^{2}}{iL+1}}~dx\\
& \qquad\quad\int_{\mathbb{R}}e^{it\tau}\left(  \frac{1}{\sqrt{-ihL+1}%
}~e^{-\frac{1}{4}\frac{\left(  t+2\tau hL\right)  ^{2}}{-ihL+1}}\right)
\left[  \left(  \frac{b_{j}^{0}}{2}+\frac{b_{j}^{1}}{2i\sqrt{\mu_{j}}}\right)
e^{it\sqrt{\mu_{j}}}\right. \\
\multicolumn{1}{r}{} & \multicolumn{1}{r}{\left.  +\left(  \frac{b_{j}^{0}}%
{2}-\frac{b_{j}^{1}}{2i\sqrt{\mu_{j}}}\right)  e^{-it\sqrt{\mu_{j}}}\right]
dt~d\tau\text{ .}}%
\end{array}
\tag{2.13}\label{2.13}%
\end{equation}

\noindent On another hand,
\begin{equation}%
\begin{array}
[c]{ll}%
\int_{\mathbb{R}}e^{i\left(  \tau\pm\sqrt{\mu_{j}}\right)  t}~e^{-\frac{1}%
{4}\frac{\left(  t+2\tau hL\right)  ^{2}}{-ihL+1}}~dt & =\int_{\mathbb{R}%
}e^{i\left(  \tau\pm\sqrt{\mu_{j}}\right)  \left(  t-2\tau hL\right)
}~e^{-\frac{1}{4}\frac{t^{2}}{-ihL+1}}~dt\\
& =e^{i\left(  \tau\pm\sqrt{\mu_{j}}\right)  \left(  -2\tau hL\right)  }%
\int_{\mathbb{R}}e^{i\left(  \tau\pm\sqrt{\mu_{j}}\right)  t}~e^{-\frac{1}%
{4}\frac{t^{2}}{-ihL+1}}~dt\\
& =e^{i\left(  \tau\pm\sqrt{\mu_{j}}\right)  \left(  -2\tau hL\right)
}~\widehat{e^{-\frac{1}{4}\frac{t^{2}}{-ihL+1}}}\left(  -\left(  \tau\pm
\sqrt{\mu_{j}}\right)  \right) \\
& =e^{i\left(  \tau\pm\sqrt{\mu_{j}}\right)  \left(  -2\tau hL\right)
}~2\sqrt{\pi}\sqrt{-ihL+1}~e^{-\left(  -ihL+1\right)  \left(  \tau\pm\sqrt
{\mu_{j}}\right)  ^{2}}%
\end{array}
\tag{2.14}\label{2.14}%
\end{equation}

\noindent where we have used the following formula
\[
\widehat{e^{-\frac{z}{2}t^{2}}}\left(  \tau\right)  =\int_{\mathbb{R}%
}e^{-it\tau}~e^{-\frac{z}{2}t^{2}}dt=\frac{\sqrt{2\pi}}{\sqrt{z}}~e^{-\frac
{1}{2z}\tau^{2}}\text{,\quad}\operatorname{Re}z\geq0\text{ and }z\neq0\text{
.}%
\]

\noindent Consequently, by (\ref{2.13})-(\ref{2.14}),
\[%
\begin{array}
[c]{ll}%
& \quad\int_{\Omega\times\mathbb{R}}A_{L}^{n}\left(  x_{o},\xi_{o3}\right)
f\left(  x,t\right)  u\left(  x,t\right)  dxdt\\
& =\frac{1}{\left(  iL+1\right)  ^{3/2}}\frac{\left(  -1\right)  ^{n}}{\left(
2\pi\right)  ^{4}}\sum\limits_{j\geq1}\int_{\Omega}\ell_{j}\left(  x\right) \\
& \qquad\int_{\mathbb{R}^{2}}\int_{\xi_{o3}-1}^{\xi_{o3}+1}\int_{\left|
\tau\right|  <\lambda}e^{i\left(  x_{1}\xi_{1}+x_{2}\xi_{2}\right)
}~e^{i\left[  \left(  -1\right)  ^{n}x_{3}+2n\frac{\xi_{o3}}{\left|  \xi
_{o3}\right|  }\rho\right]  \xi_{3}}~e^{-i\left(  \xi^{2}-\tau^{2}\right)
hL}~\widehat{\varphi f}\left(  \xi,\tau\right) \\
& \qquad\quad e^{-\frac{1}{4h}\frac{\left(  x_{1}-x_{o1}-2\xi_{1}hL\right)
^{2}}{iL+1}}~e^{-\frac{1}{4h}\frac{\left(  x_{2}-x_{o2}-2\xi_{2}hL\right)
^{2}}{iL+1}}~e^{-\frac{1}{4h}\frac{\left(  \left(  -1\right)  ^{n}%
x_{3}+2n\frac{\xi_{o3}}{\left|  \xi_{o3}\right|  }\rho-x_{o3}-2\xi
_{3}hL\right)  ^{2}}{iL+1}}~dxd\xi\\
& \qquad\quad\left[  \left(  \frac{b_{j}^{0}}{2}+\frac{b_{j}^{1}}{2i\sqrt
{\mu_{j}}}\right)  e^{i\left(  \tau+\sqrt{\mu_{j}}\right)  \left(  -2\tau
hL\right)  }~2\sqrt{\pi}e^{-\left(  -ihL+1\right)  \left(  \tau+\sqrt{\mu_{j}%
}\right)  ^{2}}\right. \\
& \qquad\quad\quad+\left.  \left(  \frac{b_{j}^{0}}{2}-\frac{b_{j}^{1}%
}{2i\sqrt{\mu_{j}}}\right)  e^{i\left(  \tau-\sqrt{\mu_{j}}\right)  \left(
-2\tau hL\right)  }~2\sqrt{\pi}e^{-\left(  -ihL+1\right)  \left(  \tau
-\sqrt{\mu_{j}}\right)  ^{2}}\right]  ~d\tau\text{ .}%
\end{array}
\]

\noindent Thus,
\[%
\begin{array}
[c]{ll}%
& \quad\left|  \int_{\Omega\times\mathbb{R}}A_{L}^{n}\left(  x_{o},\xi
_{o}\right)  f\left(  x,t\right)  u\left(  x,t\right)  dxdt\right| \\
& \leq\frac{1}{\left(  \sqrt{L^{2}+1}\right)  ^{3/2}}\frac{2\sqrt{\pi}%
}{\left(  2\pi\right)  ^{4}}\sum\limits_{j\geq1}\int_{\Omega}\left|  \ell
_{j}\left(  x\right)  \right|  \int_{\mathbb{R}^{2}}\int_{\xi_{o3}-1}%
^{\xi_{o3}+1}\int_{\left|  \tau\right|  <\lambda}\left|  \widehat{\varphi
f}\left(  \xi,\tau\right)  \right| \\
& \qquad\quad e^{-\frac{1}{4h}\left(  \left(  -1\right)  ^{n}x_{3}+2n\frac
{\xi_{o3}}{\left|  \xi_{o3}\right|  }\rho-x_{o3}-2\xi_{3}hL\right)  ^{2}%
\frac{1}{L^{2}+1}}~dx\\
& \qquad\quad\left(  \left|  \frac{b_{j}^{0}}{2}+\frac{b_{j}^{1}}{2i\sqrt
{\mu_{j}}}\right|  +\left|  \frac{b_{j}^{0}}{2}-\frac{b_{j}^{1}}{2i\sqrt
{\mu_{j}}}\right|  \right)  d\xi d\tau\text{ .}%
\end{array}
\]

\noindent Then,
\begin{equation}%
\begin{array}
[c]{ll}%
& \quad\sum\limits_{n\in\mathbb{Z}}\left|  \int_{\Omega\times\mathbb{R}}%
A_{L}^{n}\left(  x_{o},\xi_{o3}\right)  f\left(  x,t\right)  u\left(
x,t\right)  dxdt\right| \\
& \leq c\frac{1}{\left(  \sqrt{L^{2}+1}\right)  ^{3/2}}\left(  4+c\sqrt
{h\left(  L^{2}+1\right)  }\right)  \sum\limits_{j\geq1}\left(  \left|
b_{j}^{0}\right|  +\left|  \frac{b_{j}^{1}}{\sqrt{\mu_{j}}}\right|  \right)
\int_{\mathbb{R}^{2}}\int_{\xi_{o3}-1}^{\xi_{o3}+1}\int_{\left|  \tau\right|
<\lambda}\left|  \widehat{\varphi f}\left(  \xi,\tau\right)  \right|  d\xi
d\tau
\end{array}
\tag{2.15}\label{2.15}%
\end{equation}

\noindent by using
\[
\sum\limits_{n\in\mathbb{Z}}e^{-\frac{1}{4h}\left(  \left(  -1\right)
^{n}x_{3}\frac{\left|  \xi_{o3}\right|  }{\xi_{o3}}+2n\rho-x_{o3}\frac{\left|
\xi_{o3}\right|  }{\xi_{o3}}-2\xi_{3}hL\frac{\left|  \xi_{o3}\right|  }%
{\xi_{o3}}\right)  ^{2}\frac{1}{L^{2}+1}}\leq4+c\sqrt{h\left(  L^{2}+1\right)
}\text{ ,}%
\]

\noindent and $\int_{\Omega}\left|  \ell_{j}\left(  x\right)  \right|  dx\leq
c\left\|  \ell_{j}\right\|  _{L^{2}\left(  \Omega\right)  }=c$. On another
hand, using Cauchy-Schwartz inequality, we have
\begin{equation}%
\begin{array}
[c]{ll}%
\sum\limits_{j\geq1}\left(  \left|  b_{j}^{0}\right|  +\left|  \frac{b_{j}%
^{1}}{\sqrt{\mu_{j}}}\right|  \right)  & \leq\sqrt{\sum\limits_{j\geq1}\mu
_{j}^{2}\left(  \left|  b_{j}^{0}\right|  +\left|  \frac{b_{j}^{1}}{\sqrt
{\mu_{j}}}\right|  \right)  ^{2}}\sqrt{\sum\limits_{j\geq1}\left|  \frac
{1}{\mu_{j}}\right|  ^{2}}\\
& \leq c\sqrt{\mathcal{G}\left(  \partial_{t}u,0\right)  }%
\end{array}
\tag{2.16}\label{2.16}%
\end{equation}

\noindent because
\[
\sum\limits_{j\geq1}\left|  \frac{1}{\mu_{j}}\right|  ^{2}\leq c\sum
\limits_{j\geq1}\left|  \frac{1}{j^{2/3}}\right|  ^{2}<+\infty\text{ .}%
\]

\noindent We conclude from (\ref{2.12}), (\ref{2.15})-(\ref{2.16}) and
(\ref{2.8}), that for any $\left(  P,Q\right)  \in\mathbb{N}^{2}$,
\begin{equation}%
\begin{array}
[c]{ll}%
& \quad\int_{\Omega\times\mathbb{R}}\sum\limits_{\xi_{o3}\in\left(
2\mathbb{Z}+1\right)  }A_{L,P,Q}\left(  x_{o},\xi_{o3}\right)  f\left(
x,t\right)  u\left(  x,t\right)  dxdt\\
& \leq c\frac{1}{1+L\sqrt{L}}\left(  1+\sqrt{h}L\right)  \int_{\mathbb{R}^{3}%
}\int_{\left|  \tau\right|  <\lambda}\left|  \widehat{\varphi f}\left(
\xi,\tau\right)  \right|  d\xi d\tau~\sqrt{\mathcal{G}\left(  \partial
_{t}u,0\right)  }\\
& \leq c\frac{1}{\sqrt{L}}\left(  \frac{\lambda}{h}\right)  ^{\gamma}%
\sqrt{\mathcal{G}\left(  u,0\right)  }\sqrt{\mathcal{G}\left(  \partial
_{t}u,0\right)  }\text{ .}%
\end{array}
\tag{2.17}\label{2.17}%
\end{equation}

\bigskip

\subsection{The boundary term}

\bigskip

In this subsection, we study the boundary term appearing in (\ref{2.11})%
\[
ih\int_{0}^{L}\int_{\partial\Omega\times\mathbb{R}}\sum\limits_{\xi_{o3}%
\in\left(  2\mathbb{Z}+1\right)  }A_{s,P\left(  \xi_{o3}\right)  ,Q\left(
\xi_{o3}\right)  }\left(  x_{o},\xi_{o3}\right)  f\left(  x,t\right)
\partial_{\nu}u\left(  x,t\right)  dxdtds\text{ .}%
\]

\noindent Recall that $\partial\Omega=\Gamma_{1}\cup\Gamma_{2}\cup\Upsilon$.
We begin to estimate
\[
ih\int_{0}^{L}\int_{\left(  \Gamma_{1}\cup\Gamma_{2}\right)  \times\mathbb{R}%
}\sum\limits_{\xi_{o3}\in\left(  2\mathbb{Z}+1\right)  }A_{s,P\left(  \xi
_{o3}\right)  ,Q\left(  \xi_{o3}\right)  }\left(  x_{o},\xi_{o3}\right)
f\left(  x,t\right)  \partial_{\nu}u\left(  x,t\right)  dxdtds\text{ .}%
\]

\noindent First, it holds
\begin{equation}%
\begin{array}
[c]{ll}%
& \quad\int_{0}^{L}\int_{\left(  \Gamma_{1}\cup\Gamma_{2}\right)
\times\mathbb{R}}\sum\limits_{\xi_{o3}\in\left(  2\mathbb{Z}+1\right)
}A_{s,P\left(  \xi_{o3}\right)  ,Q\left(  \xi_{o3}\right)  }\left(  x_{o}%
,\xi_{o3}\right)  f\left(  x,t\right)  \partial_{\nu}u\left(  x,t\right)
dxdtds\\
& \leq\int_{0}^{L}\int_{-m_{1}}^{m_{1}}\int_{-m_{2}}^{m_{2}}\int_{\mathbb{R}%
}\sum\limits_{\xi_{o3}\in\left(  2\mathbb{Z}+1\right)  }\left|  A_{s,P\left(
\xi_{o3}\right)  ,Q\left(  \xi_{o3}\right)  }\left(  x_{o},\xi_{o3}\right)
f\left(  x_{1},x_{2},-\frac{\xi_{o3}}{\left|  \xi_{o3}\right|  }\rho,t\right)
\right| \\
& \qquad\qquad\qquad\qquad\qquad\qquad\qquad\qquad\qquad\qquad\qquad\left|
\partial_{x_{3}}u\left(  x_{1},x_{2},-\frac{\xi_{o3}}{\left|  \xi_{o3}\right|
}\rho,t\right)  \right|  dx_{1}dx_{2}dtds\\
& \quad+\int_{0}^{L}\int_{-m_{1}}^{m_{1}}\int_{-m_{2}}^{m_{2}}\int
_{\mathbb{R}}\sum\limits_{\xi_{o3}\in\left(  2\mathbb{Z}+1\right)  }\left|
A_{s,P\left(  \xi_{o3}\right)  ,Q\left(  \xi_{o3}\right)  }\left(  x_{o}%
,\xi_{o3}\right)  f\left(  x_{1},x_{2},\frac{\xi_{o3}}{\left|  \xi
_{o3}\right|  }\rho,t\right)  \right| \\
& \qquad\qquad\qquad\qquad\qquad\qquad\qquad\qquad\qquad\qquad\qquad\left|
\partial_{x_{3}}u\left(  x_{1},x_{2},\frac{\xi_{o3}}{\left|  \xi_{o3}\right|
}\rho,t\right)  \right|  dx_{1}dx_{2}dtds\text{ .}%
\end{array}
\tag{2.18}\label{2.18}%
\end{equation}

\noindent Recall that for any $\left(  P,Q\right)  \in\mathbb{N}^{2}$,
\begin{equation}%
\begin{array}
[c]{ll}%
& \quad A_{s,P,Q}\left(  x_{o},\xi_{o3}\right)  f\left(  x_{1},x_{2},\frac
{\xi_{o3}}{\left|  \xi_{o3}\right|  }\rho,t\right)  =0\quad\forall
s\geq0\text{ ,}%
\end{array}
\tag{2.19}\label{2.19}%
\end{equation}%

\[%
\begin{array}
[c]{ll}%
& \quad A_{s,P,Q}\left(  x_{o},\xi_{o3}\right)  f\left(  x_{1},x_{2}%
,-\frac{\xi_{o3}}{\left|  \xi_{o3}\right|  }\rho,t\right) \\
& =A_{s}^{-2Q}\left(  x_{o},\xi_{o3}\right)  f\left(  x_{1},x_{2},-\frac
{\xi_{o3}}{\left|  \xi_{o3}\right|  }\rho,t\right)  +A_{s}^{2P+1}\left(
x_{o},\xi_{o3}\right)  f\left(  x_{1},x_{2},-\frac{\xi_{o3}}{\left|  \xi
_{o3}\right|  }\rho,t\right)  \quad\forall s\geq0\text{ ,}%
\end{array}
\]

\noindent where%

\[%
\begin{array}
[c]{ll}%
& A_{s}^{-2Q}\left(  x_{o},\xi_{o3}\right)  f\left(  x_{1},x_{2},-\frac
{\xi_{o3}}{\left|  \xi_{o3}\right|  }\rho,t\right) \\
& \leq\frac{1}{\left(  2\pi\right)  ^{4}}\int_{\mathbb{R}^{2}}\int_{\xi
_{o3}-1}^{\xi_{o3}+1}\int_{\left|  \tau\right|  <\lambda}\left|
\widehat{\varphi f}\left(  \xi,\tau\right)  \right|  ~\left|  a\left(
0,0,\left(  4Q+1\right)  \frac{\xi_{o3}}{\left|  \xi_{o3}\right|  }\rho
+x_{o3}+2\xi_{3}hs,t+2\tau hs,s\right)  \right|  ~d\xi d\tau\text{ ,}%
\end{array}
\]%

\[%
\begin{array}
[c]{ll}%
& A_{s}^{2P+1}\left(  x_{o},\xi_{o3}\right)  f\left(  x_{1},x_{2},-\frac
{\xi_{o3}}{\left|  \xi_{o3}\right|  }\rho,t\right) \\
& \leq\frac{1}{\left(  2\pi\right)  ^{4}}\int_{\mathbb{R}^{2}}\int_{\xi
_{o3}-1}^{\xi_{o3}+1}\int_{\left|  \tau\right|  <\lambda}\left|
\widehat{\varphi f}\left(  \xi,\tau\right)  \right|  ~\left|  a\left(
0,0,\left(  4P+3\right)  \frac{\xi_{o3}}{\left|  \xi_{o3}\right|  }\rho
-x_{o3}-2\xi_{3}hs,t+2\tau hs,s\right)  \right|  ~d\xi d\tau\text{ .}%
\end{array}
\]

\noindent Now, for $L>0$, $\xi_{o3}\in\left(  2\mathbb{Z}+1\right)  $ and
$\left(  x_{o3},\xi_{3}\right)  \in\left[  -\frac{\rho}{4},\frac{\rho}%
{4}\right]  \times\left[  \xi_{o3}-1,\xi_{o3}+1\right]  $, we choose $Q\left(
\xi_{o3}\right)  $ and $P\left(  \xi_{o3}\right)  $ large enough, for example%
\[
\left\{
\begin{array}
[c]{ll}%
Q\left(  \xi_{o3}\right)  = & Q=\frac{1}{4\rho}\left(  L+1\right) \\
P\left(  \xi_{o3}\right)  = & P=\frac{1}{4\rho}\left[  \left(  L+1\right)
+2\left(  \left|  \xi_{o3}\right|  +1\right)  h_{o}L\right]
\end{array}
\right.
\]

\noindent in order that for any $s\in\left[  0,L\right]  $, $h\in\left(
0,h_{o}\right]  $ and $\left(  x_{o3},\xi_{3}\right)  \in\left[  -\frac{\rho
}{4},\frac{\rho}{4}\right]  \times\left[  \xi_{o3}-1,\xi_{o3}+1\right]  $,%
\[
\left\{
\begin{array}
[c]{rr}%
& \left|  a\left(  0,0,\left(  4Q+1\right)  \frac{\xi_{o3}}{\left|  \xi
_{o3}\right|  }\rho+x_{o3}+2\xi_{3}hs,0,s\right)  \right|  \leq\frac
{1}{\left(  \sqrt{s^{2}+1}\right)  ^{3/2}}\frac{1}{\left(  \sqrt{\left(
hs\right)  ^{2}+1}\right)  ^{1/2}}~e^{-\frac{1}{4h}}\\
& \left|  a\left(  0,0,-\left(  4P+3\right)  \frac{\xi_{o3}}{\left|  \xi
_{o3}\right|  }\rho+x_{o3}+2\xi_{3}hs,0,s\right)  \right|  \leq\frac
{1}{\left(  \sqrt{s^{2}+1}\right)  ^{3/2}}\frac{1}{\left(  \sqrt{\left(
hs\right)  ^{2}+1}\right)  ^{1/2}}~e^{-\frac{1}{4h}}\text{ .}%
\end{array}
\right.
\]

\noindent Consequently,
\[%
\begin{array}
[c]{ll}%
& \quad\sum\limits_{\xi_{o3}\in\left(  2\mathbb{Z}+1\right)  }\left|
A_{s,P\left(  \xi_{o3}\right)  ,Q\left(  \xi_{o3}\right)  }\left(  x_{o}%
,\xi_{o3}\right)  f\left(  x_{1},x_{2},-\frac{\xi_{o3}}{\left|  \xi
_{o3}\right|  }\rho,t\right)  \right|  \left|  \partial_{x_{3}}u\left(
x_{1},x_{2},-\frac{\xi_{o3}}{\left|  \xi_{o3}\right|  }\rho,t\right)  \right|
\\
& \leq e^{-\frac{1}{4h}}~\frac{1}{\left(  \sqrt{s^{2}+1}\right)  ^{3/2}}%
\frac{1}{\left(  \sqrt{\left(  hs\right)  ^{2}+1}\right)  ^{1/2}}%
~\int_{\mathbb{R}^{3}}\int_{\left|  \tau\right|  <\lambda}\left|
\widehat{\varphi f}\left(  \xi,\tau\right)  \right|  ~e^{-\frac{\left(
t+2\tau hs\right)  ^{2}}{4}\frac{1}{\left(  hs\right)  ^{2}+1}}~d\xi d\tau\\
& \quad\quad\left(  \left|  \partial_{x_{3}}u\left(  x_{1},x_{2}%
,-\rho,t\right)  \right|  +\left|  \partial_{x_{3}}u\left(  x_{1},x_{2}%
,\rho,t\right)  \right|  \right)  \text{ ,}%
\end{array}
\]

\noindent and we conclude, using (\ref{2.18})-(\ref{2.19}), a classical trace
theorem and (\ref{2.8}), that%
\begin{equation}%
\begin{array}
[c]{ll}%
& \quad\left|  ih\int_{0}^{L}\int_{\left(  \Gamma_{1}\cup\Gamma_{2}\right)
\times\mathbb{R}}\sum\limits_{\xi_{o3}\in\left(  2\mathbb{Z}+1\right)
}A_{s,P\left(  \xi_{o3}\right)  ,Q\left(  \xi_{o3}\right)  }\left(  x_{o}%
,\xi_{o3}\right)  f\left(  x,t\right)  \partial_{\nu}u\left(  x,t\right)
dxdtds\right| \\
& \leq e^{-\frac{1}{4h}}~c\int_{0}^{L}\frac{1}{1+s\sqrt{s}}\frac{1}{1+\sqrt
{h}\sqrt{s}}\int_{\mathbb{R}}\int_{\mathbb{R}^{3}}\int_{\left|  \tau\right|
<\lambda}\left|  \widehat{\varphi f}\left(  \xi,\tau\right)  \right|
~e^{-\frac{\left(  t+2\tau hs\right)  ^{2}}{4}\frac{1}{\left(  hs\right)
^{2}+1}}~d\xi d\tau~\left\|  \partial_{\nu}u\left(  \cdot,t\right)  \right\|
_{L^{2}\left(  \partial\Omega\right)  }dtds\\
& \leq e^{-\frac{1}{4h}}~c\int_{0}^{L}\frac{1}{1+s\sqrt{s}}\frac{1}{1+\sqrt
{h}\sqrt{s}}\int_{\mathbb{R}}\int_{\mathbb{R}^{3}}\int_{\left|  \tau\right|
<\lambda}\left|  \widehat{\varphi f}\left(  \xi,\tau\right)  \right|
~e^{-\frac{\left(  t+2\tau hs\right)  ^{2}}{4}\frac{1}{\left(  hs\right)
^{2}+1}}~d\xi d\tau\sqrt{\mathcal{G}\left(  \partial_{t}u,t\right)  }~dtds\\
& \leq e^{-\frac{1}{4h}}~c\sqrt{\mathcal{G}\left(  \partial_{t}u,0\right)
}\int_{0}^{L}\frac{1}{1+s\sqrt{s}}\frac{1}{1+\sqrt{h}\sqrt{s}}\left(
\int_{\mathbb{R}}e^{-\frac{t^{2}}{4}\frac{1}{\left(  hs\right)  ^{2}+1}%
}dt\right)  ds~\int_{\mathbb{R}^{3}}\int_{\left|  \tau\right|  <\lambda
}\left|  \widehat{\varphi f}\left(  \xi,\tau\right)  \right|  ~d\xi d\tau\\
& \leq e^{-\frac{1}{4h}}~L~c\sqrt{\mathcal{G}\left(  \partial_{t}u,0\right)
}~\left(  \frac{\lambda}{h}\right)  ^{\gamma}\sqrt{\mathcal{G}\left(
u,0\right)  }\text{ .}%
\end{array}
\tag{2.20}\label{2.20}%
\end{equation}

\bigskip

Now, we study to following boundary term%
\[
ih\int_{0}^{L}\int_{\Upsilon\times\mathbb{R}}\sum\limits_{\xi_{o3}\in\left(
2\mathbb{Z}+1\right)  }A_{s,P\left(  \xi_{o3}\right)  ,Q\left(  \xi
_{o3}\right)  }\left(  x_{o},\xi_{o3}\right)  f\left(  x,t\right)
\partial_{\nu}u\left(  x,t\right)  dxdtds\text{ .}%
\]

\noindent First, it holds
\begin{equation}
A_{s,P,Q}\left(  x_{o},\xi_{o3}\right)  f\left(  x,t\right)  \leq
\sum\limits_{n\in\mathbb{Z}}\left|  A_{s}^{n}\left(  x_{o},\xi_{o3}\right)
f\left(  x,t\right)  \right|  \text{ .} \tag{2.21}\label{2.21}%
\end{equation}

\noindent For $\xi_{o3}\in\left(  2\mathbb{Z}+1\right)  $, we have
\[%
\begin{array}
[c]{ll}%
& \quad A_{s}^{n}\left(  x_{o},\xi_{o3}\right)  f\left(  x,t\right) \\
& \leq\frac{1}{\left(  2\pi\right)  ^{4}}\frac{1}{\left(  \sqrt{s^{2}%
+1}\right)  ^{3/2}}\frac{1}{\left(  \sqrt{\left(  hs\right)  ^{2}+1}\right)
^{1/2}}\int_{\mathbb{R}^{2}}\int_{\xi_{o3}-1}^{\xi_{o3}+1}\int_{\left|
\tau\right|  <\lambda}\left|  \widehat{\varphi f}\left(  \xi,\tau\right)
\right| \\
& \qquad e^{-\frac{1}{4h}\left(  \left(  -1\right)  ^{n}x_{3}\frac{\left|
\xi_{o3}\right|  }{\xi_{o3}}+2n\rho-x_{o3}\frac{\left|  \xi_{o3}\right|  }%
{\xi_{o3}}-2\xi_{3}hs\frac{\left|  \xi_{o3}\right|  }{\xi_{o3}}\right)
^{2}\frac{1}{s^{2}+1}}~e^{-\frac{\left(  t+2\tau hs\right)  ^{2}}{4}\frac
{1}{\left(  hs\right)  ^{2}+1}}~d\xi d\tau
\end{array}
\]

\noindent Noticing
\[
\sum\limits_{n\in\mathbb{Z}}e^{-\frac{1}{4h}\left(  \left(  -1\right)
^{n}x_{3}\frac{\left|  \xi_{o3}\right|  }{\xi_{o3}}+2n\rho-x_{o3}\frac{\left|
\xi_{o3}\right|  }{\xi_{o3}}-2\xi_{3}hs\frac{\left|  \xi_{o3}\right|  }%
{\xi_{o3}}\right)  ^{2}\frac{1}{s^{2}+1}}\leq4+c\sqrt{h\left(  s^{2}+1\right)
}\text{ ,}%
\]

\noindent we deduce that
\begin{equation}%
\begin{array}
[c]{ll}%
& \quad\sum\limits_{n\in\mathbb{Z}}\left|  A_{s}^{n}\left(  x_{o},\xi
_{o3}\right)  f\left(  x,t\right)  \right| \\
& \leq c\left(  1+\sqrt{h}s\right)  \frac{1}{1+s\sqrt{s}}\frac{1}{1+\sqrt
{h}\sqrt{s}}\int_{\mathbb{R}^{2}}\int_{\xi_{o3}-1}^{\xi_{o3}+1}\int_{\left|
\tau\right|  <\lambda}\left|  \widehat{\varphi f}\left(  \xi,\tau\right)
\right|  e^{-\frac{\left(  t+2\tau hs\right)  ^{2}}{4}\frac{1}{\left(
hs\right)  ^{2}+1}}~d\xi d\tau\text{ .}%
\end{array}
\tag{2.22}\label{2.22}%
\end{equation}

\noindent On another hand, we get%
\begin{equation}%
\begin{array}
[c]{ll}%
& \quad\int_{\Upsilon\times\mathbb{R}}e^{-\frac{\left(  t+2\tau hs\right)
^{2}}{4}\frac{1}{\left(  hs\right)  ^{2}+1}}~\left|  \partial_{\nu}u\left(
x,t\right)  \right|  dxdt\\
& \leq\int_{\Upsilon}\int_{\left|  t+2\tau hs\right|  \leq\sqrt{\frac{\left(
hs\right)  ^{2}+1}{h}}}\left|  \partial_{\nu}u\left(  x,t\right)  \right|
dxdt+e^{-\frac{1}{8h}}\int_{\Upsilon\times\mathbb{R}}e^{-\frac{\left(  t+2\tau
hs\right)  ^{2}}{8}\frac{1}{\left(  hs\right)  ^{2}+1}}~\left|  \partial_{\nu
}u\left(  x,t\right)  \right|  dxdt\\
& \leq\int_{\Upsilon}\int_{\left|  t\right|  \leq\left|  2\tau hs\right|
+\sqrt{h}s+\frac{1}{\sqrt{h}}}\left|  \partial_{\nu}u\left(  x,t\right)
\right|  dxdt+c\left(  1+hs\right)  e^{-\frac{1}{8h}}\sqrt{\mathcal{G}\left(
\partial_{t}u,0\right)  }%
\end{array}
\tag{2.23}\label{2.23}%
\end{equation}

\noindent by cutting the integral over $t\in\mathbb{R}$ into two parts and
using a classical trace theorem. From (\ref{2.21}), (\ref{2.22}), (\ref{2.23})
and (\ref{2.8}), we conclude that%
\begin{equation}%
\begin{array}
[c]{ll}%
& \quad\left|  ih\int_{0}^{L}\int_{\Upsilon\times\mathbb{R}}\sum
\limits_{\xi_{o3}\in\left(  2\mathbb{Z}+1\right)  }A_{s,P,Q}\left(  x_{o}%
,\xi_{o3}\right)  f\left(  x,t\right)  \partial_{\nu}u\left(  x,t\right)
dxdtds\right| \\
& \leq h\int_{0}^{L}\int_{\Upsilon\times\mathbb{R}}\sum\limits_{\xi_{o3}%
\in\left(  2\mathbb{Z}+1\right)  }\sum\limits_{n\in\mathbb{Z}}\left|
A_{s}^{n}\left(  x_{o},\xi_{o3}\right)  f\left(  x,t\right)  \right|  \left|
\partial_{\nu}u\left(  x,t\right)  \right|  dxdtds\\
& \leq h\int_{0}^{L}c\left(  1+\sqrt{h}s\right)  \frac{1}{1+s\sqrt{s}}\frac
{1}{1+\sqrt{h}\sqrt{s}}\sum\limits_{\xi_{o3}\in\left(  2\mathbb{Z}+1\right)
}\int_{\mathbb{R}^{2}}\int_{\xi_{o3}-1}^{\xi_{o3}+1}\int_{\left|  \tau\right|
<\lambda}\left|  \widehat{\varphi f}\left(  \xi,\tau\right)  \right|  d\xi\\
& \qquad\qquad\left(  \int_{\Upsilon\times\mathbb{R}}e^{-\frac{\left(  t+2\tau
hs\right)  ^{2}}{4}\frac{1}{\left(  hs\right)  ^{2}+1}}\left|  \partial_{\nu
}u\left(  x,t\right)  \right|  dxdt\right)  ~d\tau ds\\
& \leq h\int_{0}^{L}c\left(  1+\sqrt{h}s\right)  \frac{1}{1+s\sqrt{s}}\frac
{1}{1+\sqrt{h}\sqrt{s}}\left(  \left(  \frac{\lambda}{h}\right)  ^{\gamma
}\sqrt{\mathcal{G}\left(  u,0\right)  }\right) \\
& \qquad\qquad\left(  \int_{\Upsilon}\int_{\left|  t\right|  \leq2\lambda
hL+\sqrt{h}L+\frac{1}{\sqrt{h}}}\left|  \partial_{\nu}u\left(  x,t\right)
\right|  dxdt+c\left(  1+hs\right)  ~e^{-\frac{1}{8h}}\sqrt{\mathcal{G}\left(
\partial_{t}u,0\right)  }\right)  ds\\
& \leq L~c\left(  \frac{\lambda}{h}\right)  ^{\alpha}~\left(  \int_{\Upsilon
}\int_{\left|  t\right|  \leq\alpha\left(  \frac{\lambda}{h}\right)  ^{\alpha
}}\left|  \partial_{\nu}u\left(  x,t\right)  \right|  ^{2}dxdt\right)
^{1/2}\sqrt{\mathcal{G}\left(  u,0\right)  }+L~c\left(  \frac{\lambda}%
{h}\right)  ^{\gamma}~e^{-\frac{1}{8h}}\sqrt{\mathcal{G}\left(  u,0\right)
}\sqrt{\mathcal{G}\left(  \partial_{t}u,0\right)  }\text{ .}%
\end{array}
\tag{2.24}\label{2.24}%
\end{equation}

\bigskip

\noindent Finally, (\ref{2.20}) and (\ref{2.24}) imply%
\begin{equation}%
\begin{array}
[c]{ll}%
& \quad\left|  ih\int_{0}^{L}\int_{\partial\Omega\times\mathbb{R}}%
\sum\limits_{\xi_{o3}\in\left(  2\mathbb{Z}+1\right)  }A_{s,P\left(  \xi
_{o3}\right)  ,Q\left(  \xi_{o3}\right)  }\left(  x_{o},\xi_{o3}\right)
f\left(  x,t\right)  \partial_{\nu}u\left(  x,t\right)  dxdtds\right| \\
& \leq L~c\left(  \frac{\lambda}{h}\right)  ^{\gamma}\left(  \int_{\Upsilon
}\int_{\left|  t\right|  \leq\gamma\left(  \frac{\lambda}{h}\right)  ^{\gamma
}}\left|  \partial_{\nu}u\left(  x,t\right)  \right|  ^{2}dxdt\right)
^{1/2}\sqrt{\mathcal{G}\left(  u,0\right)  }+L~c\left(  \frac{\lambda}%
{h}\right)  ^{\gamma}e^{-\frac{1}{4h}}\sqrt{\mathcal{G}\left(  u,0\right)
}\sqrt{\mathcal{G}\left(  \partial_{t}u,0\right)  }\text{ .}%
\end{array}
\tag{2.25}\label{2.25}%
\end{equation}

\bigskip

\subsection{The choice of $\lambda$ and $L$}

\bigskip

By (\ref{2.5}), (\ref{2.9}), (\ref{2.11}), (\ref{2.17}) and (\ref{2.25}), we
obtain, when $f=\partial_{t}^{2}u$,
\begin{equation}%
\begin{array}
[c]{ll}%
& \quad\int_{\Omega\times\mathbb{R}}a_{o}\left(  x,t\right)  \varphi\left(
x,t\right)  \partial_{t}^{2}u\left(  x,t\right)  u\left(  x,t\right)  dxdt\\
& \leq c\frac{1}{\sqrt{\lambda}}\sqrt{\mathcal{G}\left(  u,0\right)  }%
\sqrt{\mathcal{G}\left(  \partial_{t}u,0\right)  }\\
& \quad+c\left(  \frac{\lambda}{h}\right)  ^{\gamma}e^{-\frac{c}{h}%
}\mathcal{G}\left(  u,0\right) \\
& \quad+c\frac{1}{\sqrt{L}}~\left(  \frac{\lambda}{h}\right)  ^{\gamma}%
\sqrt{\mathcal{G}\left(  u,0\right)  }\sqrt{\mathcal{G}\left(  \partial
_{t}u,0\right)  }\\
& \quad+L~c\left(  \frac{\lambda}{h}\right)  ^{\gamma}\left(  \int_{\Upsilon
}\int_{\left|  t\right|  \leq\gamma\left(  \frac{\lambda}{h}\right)  ^{\gamma
}}\left|  \partial_{\nu}u\left(  x,t\right)  \right|  ^{2}dxdt\right)
^{1/2}\sqrt{\mathcal{G}\left(  u,0\right)  }+L~c\left(  \frac{\lambda}%
{h}\right)  ^{\gamma}e^{-\frac{c}{h}}\sqrt{\mathcal{G}\left(  u,0\right)
}\sqrt{\mathcal{G}\left(  \partial_{t}u,0\right)  }\text{ .}%
\end{array}
\tag{2.26}\label{2.26}%
\end{equation}

\noindent We choose $\lambda\geq1$ and $L\geq1$ be such that $\frac{1}%
{h\sqrt{h}}\frac{1}{\sqrt{\lambda}}=\sqrt{h}$ and $\frac{1}{h\sqrt{h}}\frac
{1}{\sqrt{L}}~\left(  \frac{\lambda}{h}\right)  ^{\gamma}=\sqrt{h}$ in order
that
\begin{equation}%
\begin{array}
[c]{ll}%
& \quad\frac{1}{h\sqrt{h}}\int_{\Omega\times\mathbb{R}}a_{o}\left(
x,t\right)  \varphi\left(  x,t\right)  \partial_{t}^{2}u\left(  x,t\right)
u\left(  x,t\right)  dxdt\\
& \leq c\sqrt{h}\sqrt{\mathcal{G}\left(  u,0\right)  }\sqrt{\mathcal{G}\left(
\partial_{t}u,0\right)  }+c\left(  \frac{1}{h}\right)  ^{\gamma}\left(
\int_{\Upsilon}\int_{\left|  t\right|  \leq\gamma\left(  \frac{1}{h}\right)
^{\gamma}}\left|  \partial_{\nu}u\left(  x,t\right)  \right|  ^{2}dxdt\right)
^{1/2}\sqrt{\mathcal{G}\left(  u,0\right)  }\text{ .}%
\end{array}
\tag{2.27}\label{2.27}%
\end{equation}

\noindent By replacing $f$ by $u$ and $a$ by $\widetilde{a}$ solution of
$\left(  i\partial_{s}+h\left(  \Delta-\partial_{t}^{2}\right)  \right)
\widetilde{a}\left(  x,t,s\right)  =0$ given by
\[
\widetilde{a}\left(  x,t,s\right)  =\left(  \frac{1}{\left(  is+1\right)
^{3/2}}~e^{-\frac{1}{4h}\frac{x^{2}}{is+1}}\right)  \left(  \frac{\sqrt{2}%
}{\sqrt{-ihs+2}}~e^{-\frac{1}{4}\frac{t^{2}}{-ihs+2}}\right)
\]
\noindent such that $\widetilde{a}\left(  x-x_{o},t,0\right)  =a_{o}\left(
x,t/\sqrt{2}\right)  $,we can argue in a similar fashion than above to show
that
\begin{equation}%
\begin{array}
[c]{ll}%
& \quad\frac{1}{h\sqrt{h}}\int_{\Omega\times\mathbb{R}}a_{o}\left(
x,t/\sqrt{2}\right)  \varphi\left(  x,t/\sqrt{2}\right)  u\left(  x,t\right)
u\left(  x,t\right)  dxdt\\
& \leq c\sqrt{h}\sqrt{\mathcal{G}\left(  u,0\right)  }\sqrt{\mathcal{G}\left(
\partial_{t}u,0\right)  }+c\left(  \frac{1}{h}\right)  ^{\gamma}\left(
\int_{\Upsilon}\int_{\left|  t\right|  \leq\gamma\left(  \frac{1}{h}\right)
^{\gamma}}\left|  \partial_{\nu}u\left(  x,t\right)  \right|  ^{2}dxdt\right)
^{1/2}\sqrt{\mathcal{G}\left(  u,0\right)  }\text{ .}%
\end{array}
\tag{2.28}\label{2.28}%
\end{equation}

\noindent Consequently, from (\ref{2.2}), (\ref{2.27})-(\ref{2.28}), we obtain
that for any $\left\{  x_{o}^{i}\right\}  _{i\in I}\in\overline{\omega_{o}}$,
\[%
\begin{array}
[c]{ll}%
& \quad\frac{1}{h\sqrt{h}}\int_{\Omega\times\mathbb{R}}\chi_{x_{o}^{i}}\left(
x\right)  \left|  a\left(  x-x_{o}^{i},t,0\right)  \partial_{t}u\left(
x,t\right)  \right|  ^{2}dxdt\\
& \leq c\sqrt{h}\sqrt{\mathcal{G}\left(  u,0\right)  }\sqrt{\mathcal{G}\left(
\partial_{t}u,0\right)  }+c\left(  \frac{1}{h}\right)  ^{\gamma}\left(
\int_{\Upsilon}\int_{\left|  t\right|  \leq\gamma\left(  \frac{1}{h}\right)
^{\gamma}}\left|  \partial_{\nu}u\left(  x,t\right)  \right|  ^{2}dxdt\right)
^{1/2}\sqrt{\mathcal{G}\left(  u,0\right)  }\text{ ,}%
\end{array}
\]

\noindent and finally, by (\ref{2.3}),
\[%
\begin{array}
[c]{ll}%
& \quad\int_{\omega_{o}\times\left(  0,T\right)  }\left|  \partial_{t}u\left(
x,t\right)  \right|  ^{2}dxdt\\
& \leq C_{o}\sqrt{h}\sqrt{\mathcal{G}\left(  u,0\right)  }\sqrt{\mathcal{G}%
\left(  \partial_{t}u,0\right)  }+C_{o}\left(  \frac{1}{h}\right)  ^{\gamma
}\left(  \int_{\Upsilon}\int_{\left|  t\right|  \leq\gamma\left(  \frac{1}%
{h}\right)  ^{\gamma}}\left|  \partial_{\nu}u\left(  x,t\right)  \right|
^{2}dxdt\right)  ^{1/2}\sqrt{\mathcal{G}\left(  u,0\right)  }\text{ .}%
\end{array}
\]

\noindent for some $C_{o}>0$ independent of $u$ and $h$.

\bigskip

\section{Proof of the main result}

\bigskip

The choice $\omega_{o}=\left(  -m_{1}+r_{o},m_{1}-r_{o}\right)  \times\left(
-m_{2}+r_{o},m_{2}-r_{o}\right)  \times\left(  -\frac{\rho}{4},\frac{\rho}%
{4}\right)  $ allows to invoke the geometrical control condition. Indeed,
observe that any generalized ray of $\partial_{t}^{2}-\Delta$ parametrized by
$t\in\left[  0,T\right]  $, for some $T<+\infty$ meets $\omega_{o}\cup\omega$.
As a result, we have the following observability estimate. There exists $C>0$
such that for any initial data $\left(  u_{0},u_{1}\right)  \in H_{0}%
^{1}\left(  \Omega\right)  \times L^{2}\left(  \Omega\right)  $, the solution
$u$ of (\ref{1.2}) satisfies%
\begin{equation}
\left\|  \left(  u_{0},u_{1}\right)  \right\|  _{H_{0}^{1}\left(
\Omega\right)  \times L^{2}\left(  \Omega\right)  }^{2}\leq C\left(
\int_{\omega}\int_{0}^{T}\left|  \partial_{t}u\left(  x,t\right)  \right|
^{2}dxdt+\int_{\omega_{o}}\int_{0}^{T}\left|  \partial_{t}u\left(  x,t\right)
\right|  ^{2}dxdt\right)  \text{ .} \tag{3.1}\label{3.1}%
\end{equation}

\noindent From Theorem 2, the solution $u$ of (\ref{1.2}) satisfies the
following interpolation inequality. There exist $C>0$ and $\gamma>1$ such that
for any $h\in\left(  0,h_{o}\right]  $,%
\begin{equation}%
\begin{array}
[c]{ll}%
\int_{\omega_{o}}\int_{0}^{T}\left|  \partial_{t}u\left(  x,t\right)  \right|
^{2}dxdt & \leq C\left(  \frac{1}{h}\right)  ^{\gamma}\left(  \int_{\Upsilon
}\int_{\left|  t\right|  \leq\gamma\left(  \frac{1}{h}\right)  ^{\gamma}%
}\left|  \partial_{\nu}u\left(  x,t\right)  \right|  ^{2}dxdt\right)
^{1/2}\left\|  \left(  u_{0},u_{1}\right)  \right\|  _{H_{0}^{1}\left(
\Omega\right)  \times L^{2}\left(  \Omega\right)  }\\
& \quad+C\sqrt{h}~\left\|  \left(  u_{0},u_{1}\right)  \right\|  _{H^{2}\cap
H_{0}^{1}\left(  \Omega\right)  \times H_{0}^{1}\left(  \Omega\right)
}\left\|  \left(  u_{0},u_{1}\right)  \right\|  _{H_{0}^{1}\left(
\Omega\right)  \times L^{2}\left(  \Omega\right)  }\text{ .}%
\end{array}
\tag{3.2}\label{3.2}%
\end{equation}

\noindent It is now known \cite{Li} using the multipliers technique that the
normal derivative of the solution of the wave equation with Dirichlet boundary
condition satisfies the following inequalities. Let $\psi\in C_{0}^{\infty
}\left(  \Theta\right)  $ be such that $\psi=1$ on $\Upsilon$, then there is
$c>0$ such that for any $T>0$,%
\begin{equation}%
\begin{array}
[c]{ll}%
\left\|  \partial_{\nu}u\right\|  _{L^{2}\left(  \Upsilon\times\left(
-T,T\right)  \right)  } & \leq c\left\|  \psi u\right\|  _{H^{1}\left(
\omega\times\left(  -2T,2T\right)  \right)  }\\
& \leq c\left(  \left\|  \partial_{t}u\right\|  _{L^{2}\left(  \omega
\times\left(  -3T,3T\right)  \right)  }+\left\|  u\right\|  _{L^{2}\left(
\omega\times\left(  -3T,3T\right)  \right)  }\right)  \text{ .}%
\end{array}
\tag{3.3}\label{3.3}%
\end{equation}

\noindent Consequently, from (\ref{3.1})-(\ref{3.2})-(\ref{3.3}), a
translation in time, we obtain that the solution $u$ of (\ref{1.2}%
)\ satisfies
\begin{equation}%
\begin{array}
[c]{ll}%
\left\|  \left(  u_{0},u_{1}\right)  \right\|  _{H_{0}^{1}\left(
\Omega\right)  \times L^{2}\left(  \Omega\right)  } & \leq C\left(  \frac
{1}{h}\right)  ^{\gamma}\left(  \int_{\omega}\int_{0}^{\gamma\left(  \frac
{1}{h}\right)  ^{\gamma}}\left|  \partial_{t}u\left(  x,t\right)  \right|
^{2}dxdt\right)  ^{1/2}\\
& \quad+C\left(  \frac{1}{h}\right)  ^{\gamma}\left(  \int_{\omega}\int
_{0}^{\gamma\left(  \frac{1}{h}\right)  ^{\gamma}}\left|  u\left(  x,t\right)
\right|  ^{2}dxdt\right)  ^{1/2}\\
& \quad+C\sqrt{h}~\left\|  \left(  u_{0},u_{1}\right)  \right\|  _{H^{2}\cap
H_{0}^{1}\left(  \Omega\right)  \times H_{0}^{1}\left(  \Omega\right)  }\text{
.}%
\end{array}
\tag{3.4}\label{3.4}%
\end{equation}

\noindent As $\partial_{t}w$, i.e., the derivative in time of $w$ solution of
(\ref{1.1}), can be seen as a solution of the wave with a second member
$-\alpha\left(  x\right)  \partial_{t}w$, (\ref{3.4}) implies with a usual
decomposition method, knowing $\alpha>0$ on $\omega$ and $\mathcal{E}\left(
w,0\right)  \leq c\mathcal{E}\left(  \partial_{t}w,0\right)  $, that there
exist $C^{\prime}>0$ and $\gamma>1$ such that
\[%
\begin{array}
[c]{ll}%
\mathcal{E}\left(  w,0\right)  +\mathcal{E}\left(  \partial_{t}w,0\right)  &
\leq C^{\prime}\left(  \frac{1}{h}\right)  ^{\gamma}\int_{\Omega}\int
_{0}^{\gamma\left(  \frac{1}{h}\right)  ^{\gamma}}\left(  \alpha\left(
x\right)  \left|  \partial_{t}^{2}w\left(  x,t\right)  \right|  ^{2}%
+\alpha\left(  x\right)  \left|  \partial_{t}w\left(  x,t\right)  \right|
^{2}\right)  dxdt\\
& \quad+C^{\prime}h~\mathcal{E}\left(  \partial_{t}^{2}w,0\right)
\text{\quad}\forall h\in\left(  0,h_{o}\right]  \text{ .}%
\end{array}
\]

\noindent This later estimate is clearly also true for any $h>h_{o}$ and some
constant $C^{\prime\prime}>0$, thus we can choose
\[
h=\frac{1}{\left(  C^{\prime}+C^{\prime\prime}\right)  }\frac{\mathcal{E}%
\left(  w,0\right)  +\mathcal{E}\left(  \partial_{t}w,0\right)  }%
{\mathcal{E}\left(  \partial_{t}^{2}w,0\right)  }%
\]

\noindent in order that there exists $C>0$ such that
\[%
\begin{array}
[c]{cc}%
\frac{\mathcal{E}\left(  w,0\right)  +\mathcal{E}\left(  \partial
_{t}w,0\right)  }{\mathcal{E}\left(  \partial_{t}^{2}w,0\right)  } & \leq
C\left(  \frac{\mathcal{E}\left(  \partial_{t}^{2}w,0\right)  }{\mathcal{E}%
\left(  w,0\right)  +\mathcal{E}\left(  \partial_{t}w,0\right)  }\right)
^{\gamma}\int_{\Omega}\int_{0}^{C\left(  \frac{\mathcal{E}\left(  \partial
_{t}^{2}w,0\right)  }{\mathcal{E}\left(  w,0\right)  +\mathcal{E}\left(
\partial_{t}w,0\right)  }\right)  ^{\gamma}}\frac{\alpha\left(  x\right)
\left|  \partial_{t}^{2}w\left(  x,t\right)  \right|  ^{2}+\alpha\left(
x\right)  \left|  \partial_{t}w\left(  x,t\right)  \right|  ^{2}}%
{\mathcal{E}\left(  \partial_{t}^{2}w,0\right)  }dxdt\text{ .}%
\end{array}
\]

\noindent By a translation in time, we obtain that
\[%
\begin{array}
[c]{ll}%
\frac{\mathcal{E}\left(  w,s\right)  +\mathcal{E}\left(  \partial
_{t}w,s\right)  }{\mathcal{E}\left(  \partial_{t}^{2}w,0\right)  } & \leq
C\left(  \frac{\mathcal{E}\left(  \partial_{t}^{2}w,0\right)  }{\mathcal{E}%
\left(  w,s\right)  +\mathcal{E}\left(  \partial_{t}w,s\right)  }\right)
^{\gamma}\int_{\Omega}\int_{s}^{s+C\left(  \frac{\mathcal{E}\left(
\partial_{t}^{2}w,0\right)  }{\mathcal{E}\left(  w,s\right)  +\mathcal{E}%
\left(  \partial_{t}w,s\right)  }\right)  ^{\gamma}}\frac{\alpha\left(
x\right)  \left|  \partial_{t}^{2}w\left(  x,t\right)  \right|  ^{2}%
+\alpha\left(  x\right)  \left|  \partial_{t}w\left(  x,t\right)  \right|
^{2}}{\mathcal{E}\left(  \partial_{t}^{2}w,0\right)  }dxdt\\
& \leq C\left(  \frac{\mathcal{E}\left(  \partial_{t}^{2}w,0\right)
}{\mathcal{E}\left(  w,s\right)  +\mathcal{E}\left(  \partial_{t}w,s\right)
}\right)  ^{\gamma}\left[  \left(  \frac{\mathcal{E}\left(  w,s\right)
+\mathcal{E}\left(  \partial_{t}w,s\right)  }{\mathcal{E}\left(  \partial
_{t}^{2}w,0\right)  }\right)  -\left(  \frac{\mathcal{E}\left(  w,t\right)
+\mathcal{E}\left(  \partial_{t}w,t\right)  }{\mathcal{E}\left(  \partial
_{t}^{2}w,0\right)  }\right)  _{\left|  t=s+C\left(  \frac{\mathcal{E}\left(
\partial_{t}^{2}w,0\right)  }{\mathcal{E}\left(  w,s\right)  +\mathcal{E}%
\left(  \partial_{t}w,s\right)  }\right)  ^{\gamma}\right.  }\right]  \text{
.}%
\end{array}
\]

\noindent Applying Lemma in Appendix B to
\[
\mathcal{F}\left(  s\right)  =\sigma~\frac{\mathcal{E}\left(  w,s\right)
+\mathcal{E}\left(  \partial_{t}w,s\right)  }{\mathcal{E}\left(  \partial
_{t}^{2}w,0\right)  }%
\]

\noindent where $\sigma>0$ is taken in order that $\mathcal{F}$ is bounded by
one, we get that there are $C>0$ and $\delta>0$, such that for any $t>0$ and
any initial data $\left(  \widetilde{w}_{0},\widetilde{w}_{1}\right)  \in
H^{3}\left(  \Omega\right)  \cap H_{0}^{1}\left(  \Omega\right)  \times
H^{2}\left(  \Omega\right)  \cap H_{0}^{1}\left(  \Omega\right)  $, the
solution $\widetilde{w}$ of
\begin{equation}
\left\{
\begin{array}
[c]{rll}%
\partial_{t}^{2}\widetilde{w}-\Delta\widetilde{w}+\alpha\left(  x\right)
\partial_{t}\widetilde{w} & =0 & \quad\text{in}~\Omega\times\mathbb{R}\\
\widetilde{w}=\Delta\widetilde{w} & =0 & \quad\text{on}~\partial\Omega
\times\mathbb{R}\\
\left(  \widetilde{w}\left(  \cdot,0\right)  \text{,}\partial_{t}\widetilde
{w}\left(  \cdot,0\right)  \right)  & =\left(  \widetilde{w}_{0},\widetilde
{w}_{1}\right)  & \quad\text{in}~\Omega\text{ ,}%
\end{array}
\right.  \tag{3.5}\label{3.5}%
\end{equation}

\noindent satisfies
\begin{equation}
\int_{\Omega}\left(  \left|  \partial_{t}^{2}\widetilde{w}\left(  x,t\right)
\right|  ^{2}+\left|  \nabla\partial_{t}\widetilde{w}\left(  x,t\right)
\right|  ^{2}\right)  dx\leq\frac{C}{t^{\delta}}\left\|  \left(  \widetilde
{w}_{0},\widetilde{w}_{1}\right)  \right\|  _{H^{3}\left(  \Omega\right)
\times H^{2}\left(  \Omega\right)  }^{2}\text{ .} \tag{3.6}\label{3.6}%
\end{equation}

\noindent Now, it is known \cite{Li} how to deduce from (\ref{3.6}) that there
are $C>0$ and $\delta>0$, such that for any $t>0$ and any initial data
$\left(  w_{0},w_{1}\right)  \in H^{2}\left(  \Omega\right)  \cap H_{0}%
^{1}\left(  \Omega\right)  \times H_{0}^{1}\left(  \Omega\right)  $, the
solution $w$ of (\ref{1.1}) satisfies,
\[
\int_{\Omega}\left(  \left|  \partial_{t}w\left(  x,t\right)  \right|
^{2}+\left|  \nabla w\left(  x,t\right)  \right|  ^{2}\right)  dx\leq\frac
{C}{t^{\delta}}\left\|  \left(  w_{0},w_{1}\right)  \right\|  _{H^{2}\left(
\Omega\right)  \times H^{1}\left(  \Omega\right)  }^{2}\text{ .}%
\]

\noindent Indeed, let $\widetilde{w}_{0}\in H^{3}\left(  \Omega\right)  \cap
H_{0}^{1}\left(  \Omega\right)  $ be such that $\Delta\widetilde{w}_{0}%
=w_{1}+\alpha w_{0}\in H_{0}^{1}\left(  \Omega\right)  $. In particular, it
comes%
\[
\left\|  \widetilde{w}_{0}\right\|  _{H^{3}\left(  \Omega\right)  }^{2}\leq
c\left\|  w_{1}+\alpha w_{0}\right\|  _{H_{0}^{1}\left(  \Omega\right)  }%
^{2}\text{ .}%
\]

\noindent Now, notice that $w$ solves
\[
\left\{
\begin{array}
[c]{c}%
\partial_{t}w\left(  x,t\right)  -w_{1}\left(  x\right)  -\int_{0}^{t}\Delta
w\left(  x,\tau\right)  d\tau+\alpha\left(  x\right)  \left(  w\left(
x,t\right)  -w_{0}\left(  x\right)  \right)  =0\text{ ,}\quad\left(
x,t\right)  \in\Omega\times\mathbb{R}_{+}\text{ ,}\\
w=0\quad\text{on}~\partial\Omega\times\mathbb{R}_{+}\text{ ,}\\
w\left(  \cdot,0\right)  =w_{0}\text{, }\partial_{t}w\left(  \cdot,0\right)
=w_{1}\quad\text{in}~\Omega\text{ ,}%
\end{array}
\right.
\]

\noindent and as a result, $\widetilde{w}\left(  x,t\right)  =\int_{0}%
^{t}w\left(  x,\tau\right)  d\tau+\widetilde{w}_{0}\left(  x\right)  $ solves
(\ref{3.5}) with $\partial_{t}\widetilde{w}\left(  \cdot,0\right)  =w_{0}$ and
satisfies (\ref{3.6}). And we conclude that%
\[
\int_{\Omega}\left(  \left|  \partial_{t}w\left(  x,t\right)  \right|
^{2}+\left|  \nabla w\left(  x,t\right)  \right|  ^{2}\right)  dx\leq\frac
{C}{t^{\delta}}\left(  c\left\|  w_{1}+\alpha w_{0}\right\|  _{H_{0}%
^{1}\left(  \Omega\right)  }^{2}+\left\|  w_{0}\right\|  _{H^{2}\left(
\Omega\right)  }^{2}\right)  \text{ .}%
\]

\noindent That completes the proof.

\bigskip

\section{Appendix A}

\bigskip

The goal of this Appendix is to prove, with the notations of the above
sections, the two following inequalities,%
\begin{equation}%
\begin{array}
[c]{ll}%
& \quad\left|  \int_{\Omega\times\mathbb{R}}a_{o}\left(  x,t\right)  \frac
{1}{\left(  2\pi\right)  ^{4}}\int_{\mathbb{R}^{3}}\int_{\left|  \tau\right|
\geq\lambda}e^{i\left(  x\xi+t\tau\right)  }~\widehat{\varphi u}\left(
\xi,\tau\right)  ~d\xi d\tau~u\left(  x,t\right)  ~dxdt\right| \\
& \quad+\left|  \int_{\Omega\times\mathbb{R}}a_{o}\left(  x,t\right)  \frac
{1}{\left(  2\pi\right)  ^{4}}\int_{\mathbb{R}^{3}}\int_{\left|  \tau\right|
\geq\lambda}e^{i\left(  x\xi+t\tau\right)  }~\widehat{\varphi\partial_{t}%
^{2}u}\left(  \xi,\tau\right)  ~d\xi d\tau~u\left(  x,t\right)  ~dxdt\right|
\\
& \leq c\sqrt{\frac{1}{\lambda}}\sqrt{\mathcal{G}\left(  u,0\right)  }%
\sqrt{\mathcal{G}\left(  \partial_{t}u,0\right)  }%
\end{array}
\tag{A1}\label{A1}%
\end{equation}

\noindent and
\begin{equation}
\int_{\mathbb{R}^{3}}\int_{\left|  \tau\right|  <\lambda}\left|
\widehat{\varphi u}\left(  \xi,\tau\right)  \right|  ~d\xi d\tau
+\int_{\mathbb{R}^{3}}\int_{\left|  \tau\right|  <\lambda}\left|
\widehat{\varphi\partial_{t}^{2}u}\left(  \xi,\tau\right)  \right|  ~d\xi
d\tau\leq c\left(  \frac{\lambda}{h}\right)  ^{\gamma}\sqrt{\mathcal{G}\left(
u,0\right)  }\text{ .} \tag{A2}\label{A2}%
\end{equation}

\bigskip

\noindent Proof of (\ref{A1}). Introduce
\[%
\begin{array}
[c]{ll}%
R\left(  g\right)  & =\left|  \int_{\Omega\times\mathbb{R}}a_{o}\left(
x,t\right)  \frac{1}{\left(  2\pi\right)  ^{4}}\int_{\mathbb{R}^{3}}%
\int_{\left|  \tau\right|  \geq\lambda}e^{i\left(  x\xi+t\tau\right)
}~\widehat{\varphi g}\left(  \xi,\tau\right)  ~d\xi d\tau~u\left(  x,t\right)
~dxdt\right| \\
& =\left|  \int_{\Omega\times\mathbb{R}}a_{o}\left(  x,t\right)  ~\partial
_{t}\left(  \frac{1}{\left(  2\pi\right)  ^{4}}\int_{\mathbb{R}^{3}}%
\int_{\left|  \tau\right|  \geq\lambda}\frac{1}{i\tau}~e^{i\left(  x\xi
+t\tau\right)  }~\widehat{\varphi g}\left(  \xi,\tau\right)  ~d\xi
d\tau\right)  ~u\left(  x,t\right)  dxdt\right| \\
& =\left|  -\int_{\Omega\times\mathbb{R}}\partial_{t}\left(  a_{o}u\left(
x,t\right)  \right)  \left(  \frac{1}{2\pi}\int_{\left|  \tau\right|
\geq\lambda}\frac{1}{i\tau}~e^{it\tau}\left[  \int_{\mathbb{R}}e^{-i\theta
\tau}\varphi g\left(  x,\theta\right)  ~d\theta\right]  d\tau\right)
~dxdt\right|  \text{ .}%
\end{array}
\]

\noindent It follows using Cauchy-Schwartz inequality and Parseval identity
that%
\[%
\begin{array}
[c]{ll}%
R\left(  g\right)  & \leq\int_{\Omega\times\mathbb{R}}\left|  \partial
_{t}\left(  a_{o}u\left(  x,t\right)  \right)  \right|  \left(  \frac{1}{2\pi
}\sqrt{\int_{\left|  \tau\right|  \geq\lambda}\frac{1}{\tau^{2}}d\tau}%
\sqrt{\int_{\mathbb{R}}\left|  \int_{\mathbb{R}}e^{-i\theta\tau}\varphi
g\left(  x,\theta\right)  d\theta\right|  ^{2}d\tau}\right)  dxdt\\
& \leq\int_{\Omega\times\mathbb{R}}\left|  \partial_{t}\left(  a_{o}u\left(
x,t\right)  \right)  \right|  \left(  \frac{1}{2\pi}\sqrt{\int_{\left|
\tau\right|  \geq\lambda}\frac{1}{\tau^{2}}d\tau}\sqrt{2\pi\int_{\mathbb{R}%
}\left|  \varphi g\left(  x,\theta\right)  \right|  ^{2}d\theta}\right)
dxdt\\
& \leq\int_{\Omega\times\mathbb{R}}\left|  \partial_{t}\left(  a_{o}u\left(
x,t\right)  \right)  \right|  \left(  \frac{1}{\sqrt{2\pi}}\sqrt{\frac
{2}{\lambda}}\left\|  \varphi g\left(  x,\cdot\right)  \right\|
_{L^{2}\left(  \mathbb{R}\right)  }\right)  dxdt\\
& \leq\frac{1}{\sqrt{\pi}}\sqrt{\frac{1}{\lambda}}~\int_{\mathbb{R}}\left\|
\partial_{t}\left(  a_{o}u\right)  \left(  \cdot,t\right)  \right\|
_{L^{2}\left(  \Omega\right)  }dt\left\|  \varphi g\left(  x,\cdot\right)
\right\|  _{L^{2}\left(  \Omega\times\mathbb{R}\right)  }\text{ .}%
\end{array}
\]

\noindent Since we have the following estimates%
\[%
\begin{array}
[c]{ll}%
\left\|  \varphi u\left(  x,\cdot\right)  \right\|  _{L^{2}\left(
\Omega\times\mathbb{R}\right)  }+\left\|  \varphi\partial_{t}^{2}u\left(
x,\cdot\right)  \right\|  _{L^{2}\left(  \Omega\times\mathbb{R}\right)  } &
\leq\sqrt{\int_{\mathbb{R}}e^{-\frac{1}{2}t^{2}}\int_{\Omega}\left|  u\left(
x,t\right)  \right|  ^{2}dxdt}\\
& \quad+\sqrt{\int_{\mathbb{R}}e^{-\frac{1}{2}t^{2}}\int_{\Omega}\left|
\partial_{t}^{2}u\left(  x,t\right)  \right|  ^{2}dxdt}\\
& \leq\sqrt{\sqrt{2\pi}}\left(  c\sqrt{\mathcal{G}\left(  u,0\right)  }%
+\sqrt{\mathcal{G}\left(  \partial_{t}u,0\right)  }\right)
\end{array}
\]%

\[%
\begin{array}
[c]{ll}%
\int_{\mathbb{R}}\left\|  \partial_{t}\left(  a_{o}u\right)  \left(
\cdot,t\right)  \right\|  _{L^{2}\left(  \Omega\right)  }dt & \leq
\int_{\mathbb{R}}\sqrt{\int_{\Omega}\left|  \partial_{t}a_{o}u\left(
x,t\right)  \right|  ^{2}dx}dt+\int_{\mathbb{R}}\sqrt{\int_{\Omega}\left|
a_{o}\partial_{t}u\left(  x,t\right)  \right|  ^{2}dx}dt\\
& \leq\int_{\mathbb{R}}\frac{\left|  t\right|  }{2}e^{-\frac{1}{4}t^{2}}%
\sqrt{\int_{\Omega}\left|  u\left(  x,t\right)  \right|  ^{2}}dxdt+\int
_{\mathbb{R}}e^{-\frac{1}{4}t^{2}}\sqrt{\int_{\Omega}\left|  \partial
_{t}u\left(  x,t\right)  \right|  ^{2}}dxdt\\
& \leq c\sqrt{\mathcal{G}\left(  u,0\right)  }\text{ , }%
\end{array}
\]

\noindent we conclude that%
\[
R\left(  u\right)  +R\left(  \partial_{t}^{2}u\right)  \leq c\sqrt{\frac
{1}{\lambda}}\sqrt{\mathcal{G}\left(  u,0\right)  }\sqrt{\mathcal{G}\left(
\partial_{t}u,0\right)  }\text{ .}%
\]

\noindent That completes the proof of (\ref{A1}).

\bigskip

Proof of (\ref{A2}). We estimate $\int_{\mathbb{R}^{3}}\int_{\left|
\tau\right|  <\lambda}\left|  \widehat{\varphi g}\left(  \xi,\tau\right)
\right|  ~d\xi d\tau$ where $g$ solves $\partial_{t}^{2}g-\Delta g=0$ in
$\Omega\times\mathbb{R}$. Recall that $\varphi\left(  \cdot,t\right)  =\chi
a_{o}\left(  \cdot,t\right)  \in C_{0}^{\infty}\left(  \Omega\right)  $. So,
using Cauchy-Schwartz inequality and Parseval identity,%
\[%
\begin{array}
[c]{ll}%
\int_{\mathbb{R}^{3}}\int_{\left|  \tau\right|  <\lambda}\left|
\widehat{\varphi g}\left(  \xi,\tau\right)  \right|  ~d\xi d\tau &
=\int_{\mathbb{R}^{3}}\int_{\left|  \tau\right|  <\lambda}\frac{1+\left|
\xi\right|  ^{2}}{1+\left|  \xi\right|  ^{2}}\left|  \widehat{\varphi
g}\left(  \xi,\tau\right)  \right|  ~d\xi d\tau\\
& \leq\int_{\left|  \tau\right|  <\lambda}\sqrt{\int_{\mathbb{R}^{3}}\frac
{1}{\left(  1+\left|  \xi\right|  ^{2}\right)  ^{2}}}d\xi\sqrt{\int
_{\mathbb{R}^{3}}\left|  \widehat{\left(  1-\Delta\right)  \left(  \varphi
g\right)  }\left(  \xi,\tau\right)  \right|  ^{2}d\xi}d\tau\\
& \leq\pi^{2}\int_{\left|  \tau\right|  <\lambda}\sqrt{\int_{\mathbb{R}^{3}%
}\left|  \widehat{\left(  1-\Delta\right)  \left(  \varphi g\right)  }\left(
\xi,\tau\right)  \right|  ^{2}d\xi}d\tau\\
& \leq\pi^{2}\sqrt{\lambda}\sqrt{\int_{\mathbb{R}^{3}}\int_{\left|
\tau\right|  <\lambda}\left|  \widehat{\left(  1-\Delta\right)  \left(
\varphi g\right)  }\left(  \xi,\tau\right)  \right|  ^{2}d\xi d\tau}\text{ .}%
\end{array}
\]

\noindent Observe that
\[%
\begin{array}
[c]{ll}%
\Delta\left(  \varphi g\right)  & =\varphi\Delta g+2\nabla\varphi\nabla
g+\Delta\varphi g\\
& =\varphi\partial_{t}^{2}g+2\varphi\left(  x_{o},t\right)  \nabla
\varphi\left(  x,0\right)  \nabla g+\Delta\varphi g\\
& =\partial_{t}^{2}\left(  \varphi g\right)  -2\partial_{t}\left(
\partial_{t}\varphi g\right)  +\left(  \Delta\varphi+\partial_{t}^{2}%
\varphi\right)  g+2\varphi\left(  x_{o},t\right)  \nabla\varphi\left(
x,0\right)  \nabla g\text{ .}%
\end{array}
\]

\noindent As a result,
\[%
\begin{array}
[c]{ll}%
\int_{\mathbb{R}^{3}}\int_{\left|  \tau\right|  <\lambda}\left|
\widehat{\varphi g}\left(  \xi,\tau\right)  \right|  ~d\xi d\tau & \leq
c\sqrt{\lambda}\lambda^{2}\sqrt{\int_{\mathbb{R}^{3}}\int_{\left|
\tau\right|  <\lambda}\sum\limits_{k=0,1,2}\left|  \widehat{\partial_{t}%
^{2-k}\varphi g}\left(  \xi,\tau\right)  \right|  ^{2}d\xi d\tau}\\
& +c\sqrt{\lambda}\sqrt{\int_{\mathbb{R}^{3}}\int_{\left|  \tau\right|
<\lambda}\left|  \widehat{\Delta\varphi g}\left(  \xi,\tau\right)  \right|
^{2}d\xi d\tau}\\
& +c\sqrt{\lambda}\sqrt{\int_{\mathbb{R}^{3}}\int_{\left|  \tau\right|
<\lambda}\left|  \widehat{\varphi\left(  x_{o},t\right)  \nabla\varphi\left(
x,0\right)  \nabla g}\left(  \xi,\tau\right)  \right|  ^{2}d\xi d\tau}\text{
.}%
\end{array}
\]

\noindent In particular, when $g=u$ we obtain, using Parseval identity%
\[%
\begin{array}
[c]{ll}%
\int_{\mathbb{R}^{3}}\int_{\left|  \tau\right|  <\lambda}\left|
\widehat{\varphi u}\left(  \xi,\tau\right)  \right|  ~d\xi d\tau & \leq
c\sqrt{\lambda}\lambda^{2}\sum\limits_{k=0,1,2}\left\|  \partial_{t}%
^{2-k}\varphi u\right\|  _{L^{2}\left(  \Omega\times\mathbb{R}\right)
}+c\sqrt{\lambda}\left\|  \Delta\varphi u\right\|  _{L^{2}\left(  \Omega
\times\mathbb{R}\right)  }\\
& \quad+c\sqrt{\lambda}\left\|  \varphi\left(  x_{o},t\right)  \nabla
\varphi\left(  x,0\right)  \nabla u\right\|  _{L^{2}\left(  \Omega
\times\mathbb{R}\right)  }\\
& \leq c\left(  \sqrt{\lambda}\lambda^{2}+c\sqrt{\lambda}\frac{1}{h}\right)
\int_{\mathbb{R}}e^{-ct^{2}}\left(  \int_{\Omega}\left|  u\left(  x,t\right)
\right|  ^{2}dx\right)  dt+c\sqrt{\lambda}\frac{1}{\sqrt{h}}\sqrt
{\mathcal{G}\left(  u,0\right)  }\text{ .}%
\end{array}
\]

\noindent On another hand, when $g=\partial_{t}^{2}u$ and using the fact that
$\left|  \tau\right|  <\lambda$,
\[%
\begin{array}
[c]{ll}%
\int_{\mathbb{R}^{3}}\int_{\left|  \tau\right|  <\lambda}\left|
\widehat{\varphi f}\left(  \xi,\tau\right)  \right|  ~d\xi d\tau & \leq
c\sqrt{\lambda}\lambda^{3}\sqrt{\int_{\mathbb{R}^{3}}\int_{\left|
\tau\right|  <\lambda}\sum\limits_{k=0,1,2,3}\left|  \widehat{\partial
_{t}^{3-k}\varphi\partial_{t}u}\left(  \xi,\tau\right)  \right|  ^{2}d\xi
d\tau}\\
& +c\sqrt{\lambda}\lambda\sqrt{\int_{\mathbb{R}^{3}}\int_{\left|  \tau\right|
<\lambda}\sum\limits_{k=0,1}\left|  \widehat{\partial_{t}^{1-k}\Delta
\varphi\partial_{t}u}\left(  \xi,\tau\right)  \right|  ^{2}d\xi d\tau}\\
& +c\sqrt{\lambda}\lambda^{2}\sqrt{\int_{\mathbb{R}^{3}}\int_{\left|
\tau\right|  <\lambda}\sum\limits_{k=0,1,2}\left|  \widehat{\partial_{t}%
^{2-k}\varphi\left(  x_{o},t\right)  \nabla\varphi\left(  x,0\right)  \nabla
u}\left(  \xi,\tau\right)  \right|  ^{2}d\xi d\tau}\\
& \leq c\left(  \sqrt{\lambda}\lambda^{3}+\sqrt{\lambda}\lambda\frac{1}%
{h}+\sqrt{\lambda}\lambda^{2}\frac{1}{\sqrt{h}}\right)  \sqrt{\mathcal{G}%
\left(  u,0\right)  }\text{ .}%
\end{array}
\]

\noindent We conclude that there exist $c>0$ and $\gamma>1$ such that for any
$h\in\left(  0,h_{o}\right]  $ and $\lambda\geq1$,
\[
\int_{\mathbb{R}^{3}}\int_{\left|  \tau\right|  <\lambda}\left|
\widehat{\varphi u}\left(  \xi,\tau\right)  \right|  ~d\xi d\tau
+\int_{\mathbb{R}^{3}}\int_{\left|  \tau\right|  <\lambda}\left|
\widehat{\varphi\partial_{t}^{2}u}\left(  \xi,\tau\right)  \right|  ~d\xi
d\tau\leq c\left(  \frac{\lambda}{h}\right)  ^{\gamma}\sqrt{\mathcal{G}\left(
u,0\right)  }\text{ .}%
\]

\noindent That completes the proof of (\ref{A2}).

\bigskip

\section{Appendix B}

\bigskip

\noindent Lemma.- \textit{Let }$\mathcal{F}$\textit{ be a continuous positive
decreasing real function on }$\left[  0,+\infty\right)  $\textit{ and bounded
by one. Suppose that there are four constants }$c_{1}>1$\textit{ and }%
$c_{2},\beta,\gamma>0$\textit{ such that }%
\[
\mathcal{F}\left(  s\right)  \leq c_{1}\left(  \frac{1}{\mathcal{F}\left(
s\right)  }\right)  ^{\beta}\left(  \mathcal{F}\left(  s\right)
-\mathcal{F}\left(  \left(  \frac{c_{2}}{\mathcal{F}\left(  s\right)
}\right)  ^{\gamma}+s\right)  \right)  \quad\forall s>0\text{ .}%
\]

\noindent\textit{Then there exist }$C>0$\textit{ and }$\delta>0$\textit{ such
that for any }$t>0$\textit{, }%
\[
\mathcal{F}\left(  t\right)  \leq\frac{C}{t^{\delta}}\text{ .}%
\]

\bigskip

\bigskip

Proof .- Let $t>0$. If $\left(  \frac{\mathcal{F}\left(  s\right)  }{c_{2}%
}\right)  ^{\gamma}<\frac{1}{t}$ then $\mathcal{F}\left(  s\right)  \leq
\frac{c_{2}}{t^{1/\gamma}}$. If $\frac{1}{t}\leq\left(  \frac{\mathcal{F}%
\left(  s\right)  }{c_{2}}\right)  ^{\gamma}$ then $\left(  \frac{c_{2}%
}{\mathcal{F}\left(  s\right)  }\right)  ^{\gamma}+s\leq t+s$, thus $\left(
\mathcal{F}\left(  t+s\right)  \right)  \leq\mathcal{F}\left(  \left(
\frac{c_{2}}{\mathcal{F}\left(  s\right)  }\right)  ^{\gamma}+s\right)  $ and
therefore%
\[
\mathcal{F}\left(  s\right)  \leq\left(  c_{1}\right)  ^{\frac{1}{\beta+1}%
}\left\{  \mathcal{F}\left(  s\right)  -\left(  \mathcal{F}\left(  t+s\right)
\right)  \right\}  ^{\frac{1}{\beta+1}}\text{ .}%
\]

\noindent Consequently,
\[
\mathcal{F}\left(  s\right)  \leq\left(  c_{1}\right)  ^{\frac{1}{\beta+1}%
}\left\{  \mathcal{F}\left(  s\right)  -\left(  \mathcal{F}\left(  t+s\right)
\right)  \right\}  ^{\frac{1}{\beta+1}}+\frac{c_{2}}{t^{1/\gamma}}\quad\forall
s,t>0\text{ .}%
\]

\noindent Let
\[
\psi_{t}\left(  s\right)  =\frac{1}{\left(  \frac{c_{1}t}{s}\right)
^{\frac{1}{\beta+1}}+\frac{c_{2}}{t^{1/\gamma}}}\text{ .}%
\]

\noindent If $\mathcal{F}\left(  s\right)  -\left(  \mathcal{F}\left(
t+s\right)  \right)  \leq\frac{t}{t+s}$, then $\mathcal{F}\left(  s\right)
\leq\left(  \frac{c_{1}t}{t+s}\right)  ^{\frac{1}{\beta+1}}+\frac{c_{2}%
}{t^{1/\gamma}}$ and thus $\psi_{t}\left(  t+s\right)  \mathcal{F}\left(
t+s\right)  \leq1$. If $\frac{t}{t+s}<\mathcal{F}\left(  s\right)  -\left(
\mathcal{F}\left(  t+s\right)  \right)  $, then $\frac{t}{t+s}\left(
\mathcal{F}\left(  s\right)  \right)  \leq\frac{t}{t+s}<\mathcal{F}\left(
s\right)  -\left(  \mathcal{F}\left(  t+s\right)  \right)  $. Therefore%
\[%
\begin{array}
[c]{ll}%
\psi_{t}\left(  t+s\right)  \mathcal{F}\left(  t+s\right)  & <\frac{s}%
{t+s}\mathcal{F}\left(  s\right)  \psi_{t}\left(  t+s\right)  =\psi_{t}\left(
s\right)  \mathcal{F}\left(  s\right)  \left(  \frac{\frac{\psi_{t}\left(
t+s\right)  }{t+s}}{\frac{\psi_{t}\left(  s\right)  }{s}}\right) \\
& <\psi_{t}\left(  s\right)  \mathcal{F}\left(  s\right)
\end{array}
\]

\noindent by using the decreasing property of $\zeta\longmapsto\frac{\psi
_{t}\left(  \zeta\right)  }{\zeta}$. We have proved that for any $s,t>0$, we
have either $\psi_{t}\left(  t+s\right)  \mathcal{F}\left(  t+s\right)  \leq
1$, or $\psi_{t}\left(  t+s\right)  \mathcal{F}\left(  t+s\right)  <\psi
_{t}\left(  s\right)  \mathcal{F}\left(  s\right)  $. In particular, we deduce
that for any $t>0$ and $n\in\mathbb{N}\left\backslash \left\{  0\right\}
\right.  $, either%
\[%
\begin{array}
[c]{ll}%
& \psi_{t}\left(  \left(  n+1\right)  t\right)  \mathcal{F}\left(  \left(
n+1\right)  t\right)  \leq1\\
\text{or} & \psi_{t}\left(  \left(  n+1\right)  t\right)  \mathcal{F}\left(
\left(  n+1\right)  t\right)  <\psi_{t}\left(  nt\right)  \mathcal{F}\left(
nt\right)  \text{ .}%
\end{array}
\]

\noindent Then inductively, it implies that
\[
\psi_{t}\left(  \left(  n+1\right)  t\right)  \mathcal{F}\left(  \left(
n+1\right)  t\right)  \leq\max\left(  1,\psi_{t}\left(  t\right)
\mathcal{F}\left(  t\right)  \right)  =1\text{ .}%
\]

\noindent Hence for all $t>0$ and $n\in\mathbb{N}\left\backslash \left\{
0\right\}  \right.  $,
\[
\mathcal{F}\left(  \left(  n+1\right)  t\right)  \leq\left(  \frac{c_{1}}%
{n+1}\right)  ^{\frac{1}{\beta+1}}+\frac{c_{2}}{t^{1/\gamma}}\text{.}%
\]

\noindent We choose $n$ such that $n+1\leq t<n+2$ and we obtain that for all
$t\geq2$,%
\[
\mathcal{F}\left(  t^{2}\right)  \leq\left(  \frac{2c_{1}}{t}\right)
^{\frac{1}{\beta+1}}+\frac{c_{2}}{t^{1/\gamma}}\text{.}%
\]

\noindent The desired result now follows immediately.

\bigskip

\bigskip

\bigskip

\bigskip
\end{document}